\documentclass[11pt]{article}   
\usepackage{amsfonts,latexsym}
\newcommand{\cleqn}{\setcounter{equation}{0}}
\newcommand{\clth}{\setcounter{theorem}{0}}
\newcommand {\sectionnew}[1]{\section{#1}\cleqn\clth}
\newcommand{\beq}{\begin{equation}}
\newcommand{\eeq}{\end{equation}}
\newcommand{\beqa}{\begin{eqnarray}}
\newcommand{\eeqa}{\end{eqnarray}}
\newcommand{\beaa}{\begin{eqnarray*}}
\newcommand{\eaa}{\end{eqnarray*}}
\newcommand{\nn}{\hfill\nonumber}
\newcommand{\text}{\textrm}
\newcommand \nc {\newcommand}
\nc \proof {\noindent {\em{Proof.\/ }}}
\nc \qed {$\Box$\hfill}
\newtheorem{theorem}{Theorem}[section]
\newtheorem{lemma}[theorem]{Lemma}
\newtheorem{proposition}[theorem]{Proposition}
\newtheorem{corollary}[theorem]{Corollary}
\newtheorem{definition}[theorem]{Definition}
\newtheorem{example}[theorem]{Example}
\newtheorem{remark}[theorem]{Remark}
\newtheorem{conjecture}[theorem]{Conjecture}
\newtheorem{question}[theorem]{Question}
\nc \bth[1] { \begin{theorem}\label{t#1} }
\nc \ble[1] { \begin{lemma}\label{l#1} }
\nc \bpr[1] { \begin{proposition}\label{p#1} }
\nc \bco[1] { \begin{corollary}\label{c#1} }
\nc \bde[1] { \begin{definition}\label{d#1}\rm }
\nc \bex[1] { \begin{example}\label{e#1}\rm }
\nc \bre[1] { \begin{remark}\label{r#1}\rm }
\nc \bcon[1] { \begin{conjecture}\label{con#1}\rm }
\nc \bque[1] { \begin{question}\label{que#1}\rm }
\nc {\eth} { \end{theorem} }
\nc {\ele} { \end{lemma} }
\nc {\epr} { \end{proposition} }
\nc {\eco} { \end{corollary} }
\nc {\ede} { \end{definition} }
\nc {\eex} { \end{example} }
\nc {\ere} { \end{remark} }
\nc {\econ} { \end{conjecture} }
\nc {\eque} { \end{question} }
\nc \eqref[1] {{\rm{(\ref{#1})}}}
\nc \thref[1]{Theorem \ref{t#1}}
\nc \leref[1]{Lemma \ref{l#1}}
\nc \prref[1]{Proposition \ref{p#1}}
\nc \coref[1]{Corollary \ref{c#1}}
\nc \deref[1]{Definition \ref{d#1}}
\nc \exref[1]{Example \ref{e#1}}
\nc \reref[1]{Remark \ref{r#1}}

\def \Wr {{\mathrm{Wr}}}
\def \a {\alpha}
\def \b {\beta}
\def \A {{\mathcal A}}

\def \d {{\mathrm d}}

\def \Cset {{\mathbb C}}
\def \Zset {{\mathbb Z}}
\def \Nset {{\mathbb N}}
\def \Vset {{\mathbb V}}
\def \ord { {\mathrm{ord}} }
\def \rank { {\mathrm{rank}} }
\def \span { {\mathrm{span}} }
\def \const { {\mathrm{const}} }
\def \mod { {\mathrm{mod}} }
\def \spec { {\mathrm{Spec}} }

\def \mult { {\mathrm{mult}} }
\def \res { {\mathrm{Res}} }

\def \Re { {\mathrm{Re}} }
\renewcommand \ker { {\mathrm{Ker}} }
\nc \Gr {Gr}
\nc \GRN { \Gr^{(N)} }
\nc \GRA[1] { \Gr_A^{(#1)} }   
\nc \GRAN { \GRA{N} }
\nc \GrA[1] { \Gr_A(#1) }
\nc \GrAa { \GrA{\alpha} }
\nc \GRB[1] { \Gr_B^{(#1)} }   
\nc \GRBN { \GRB{N} }
\nc \GrB[1] { \Gr_B(#1) }
\nc \GrBb { \GrB{\beta} }
\nc \GRMB[1] { \Gr_{MB}^{(#1)} }   
\nc \GRMBN { \GRMB{N} }
\nc \GrMB[1] { \Gr_{MB}(#1) }
\nc \GrMBb { \GrMB{\beta} }
\def\dfrac#1#2{{\displaystyle\frac{#1}{#2}}}
\begin{document}
\title{{\LARGE\bf{ Bispectral algebras of commuting ordinary differential
operators }}}
\author{
B.~Bakalov
\thanks{New address: Department of Mathematics, MIT, Cambridge, MA 02139. 
E-mail: bakalov@math.mit.edu}
\quad
E.~Horozov
\thanks{E-mail: horozov@fmi.uni-sofia.bg}
\quad
M.~Yakimov
\thanks{New address: Department of Mathematics, University of California,
Berkeley, CA 94720. 
E-mail: yakimov@math.berkeley.edu}
\\ \hfill\\ \normalsize \textit{
Department of Mathematics and Informatics, }\\
\normalsize \textit{ Sofia University, 5 J. Bourchier Blvd.,
Sofia 1126, Bulgaria }     }
\date{}
\maketitle
\begin{abstract}
We develop a systematic way for constructing bispectral algebras of 
commuting ordinary differential operators of any rank $N$. It combines and 
unifies the ideas of Duistermaat--Gr\"unbaum and Wilson.
Our construction is completely algorithmic and enables us to 
obtain all previously known classes or individual examples of bispectral 
operators. The method also provides new broad families of bispectral algebras 
which may help to penetrate deeper into the problem.
\end{abstract}
\tableofcontents
\vspace{-11cm}
\begin{flushright}
\end{flushright}
\vspace{10cm}
\setcounter{section}{-1}
\sectionnew{Introduction}
In this paper we reconsider the bispectral problem.
As stated in \cite{DG}, it asks for which ordinary
differential operators
$L(x,\partial_x)$ there exists a family of eigenfunctions $\Psi(x,z)$ 
that are
also eigenfunctions for another differential operator $\Lambda(z,\partial_z)$
but this time in the ``spectral parameter'' $z$, to wit
    \beqa
&&L(x,\partial_x) \Psi(x,z) = f(z) \Psi(x,z), \label{0.1} \hfill\\
&&\Lambda(z,\partial_z) \Psi(x,z) = \Theta(x) \Psi(x,z) \label{0.2} \hfill
     \eeqa
for some functions $f(z),\Theta(x).$ Both operators $L$ and $\Lambda$ are 
called bispectral.

This problem first appeared in \cite{G1} in connection with ``limited
angle tomography'' (see also \cite{G2, G3, DG}). Later it turned out to be
related with several, seemingly far from it, topics and in particular, with 
soliton mathematics. To be more specific, we have to mention the deep
connection with some very actively developing areas of research in 
mathematics and theoretical physics like Calogero--Moser particle system
\cite{W2, Ka} (see also \cite{R}),
additional symmetries of KdV and KP hierarchies
\cite{MZ, BHY3}, representation theory of $W_{1+\infty}$--algebra
\cite{BHY3}, etc. These studies not only
revealed the rich mathematical structure
behind the bispectral problem, but also (if we use a remark by G.~Wilson
\cite{W2}) ``deepened the mystery'' around it. Thus, not only applications,
but also purely mathematical 
questions motivated the great activity in the past few years in the bispectral 
problem.

In the present paper we construct new families of bispectral operators.
In order to explain better our contribution, we need to review some of the
achievements in the subject.

The first general result in the direction of classifying bispectral operators
belongs to J.~J.~Duistermaat and F.~A.~Gr\"unbaum \cite{DG}. They determined
all
second order operators $L$ admitting an operator $\Lambda$ such that the pair
$(L,\Lambda)$ solves the bispectral problem (\ref{0.1}, \ref{0.2}). Their
answer is as follows. If we write the operator $L$ in the standard 
Schr\"odinger form
  $$
L = \frac{\d^2}{\d x^2} + u(x),
  $$
the bispectral potentials $u(x)$ are given (up to translations and rescalings
of $x$ and $z$) by the following list:
      \beqa
&& u(x) = x \qquad\qquad\qquad\;\; {\textrm {(Airy)}};
\label{0.3}  \\
&& u(x) = c x^{-2}, \; c\in \Cset  \qquad {\textrm{(Bessel)}};
\label{0.4} \\
&& u(x), \; {\textrm{which can be obtained by finitely
many rational Darboux}} \nn\\
&&\hspace*{-1.46cm}{\textrm{transformations from}} \; u(x) = 0;
\label{0.5} \\
&& u(x), \; {\textrm {which can be obtained by finitely
many rational Darboux}}  \nn\\
&&\hspace*{-1.46cm}{\textrm{transformations from}} \; u(x) = -\frac{1}{4x^{2}}.
\label{0.6}
        \eeqa

The family \eqref{0.5} has previously appeared in \cite{AMM, AM} and is known
as ``{\em{rational solutions of KdV\/}}''. They can be obtained also
by applying ``{\em{higher KdV flows\/}}'' to potentials
$v_k(x) = k(k+1) x^{-2}, \; k \in \Nset.$

The second family \eqref{0.6} was interpreted by F.~Magri and J.~Zubelli
\cite{MZ} as potentials invariant under the flows of the
``{\em{master symmetries\/}}'' or {\em{Virasoro flows}}.

Besides the classification of the bispectral operators by their order,
another scheme has been suggested in \cite{DG} and used in \cite{W}.
Below we explain it in a general context as it will be used throughout this 
paper.
One may consider an operator $L(x,\partial_x)$ as an element of a maximal
algebra $\A$ of commuting ordinary differential operators \cite{BC}. Following
G.~Wilson \cite{W}, we call such an algebra {\it {bispectral}} 
if there exists a
joint eigenfunction $\Psi(x,z)$ for the operators $L$ in $\A$ that satisfies
also equation \eqref{0.2}. The dimension of the space of eigenfunctions
$\Psi(x,z)$ is called {\it {rank}} of the commutative algebra $\A$ (see e.g.
\cite{KrN}). This number coincides with the greatest common divisor of the
orders of the operators in $\A$. For example, the operators with potentials
\eqref{0.5} belong to rank $1$ algebras and those with potentials
(\ref{0.3}, \ref{0.4}, \ref{0.6}) to rank $2$ algebras \cite{DG}.

All rank $1$ maximal bispectral algebras were recently found by G.~Wilson
\cite{W}. These algebras do not necessarily contain an operator of order two.

The methods of the above mentioned papers \cite{DG} and \cite{W} may seem quite
different. Indeed, while in \cite{DG} the ``rational'' Darboux transformations
play a decisive role, G.~Wilson \cite{W} uses planes in Sato's Grassmannian
obtained from the standard $H_+ = \span\{z^k\}_{k \geq 0}$ by imposing a number
of conditions on it. One of our main observations is that both methods,
appropriately modified, can be looked upon as the 
two sides of one general theory.

{}From this new point of view in the
present paper we construct nontrivial maximal bispectral algebras of any rank 
$N$, thus extending the results from \cite{DG, W}.
For example, for any positive integer $k$ we obtain bispectral
algebras of rank $N$ with the lowest order of the operators equal
to $k N$. Our method allows us to obtain all 
classes and single examples of bispectral operators known to us by a unique
method. At the same time we suggest an effective procedure for constructing
bispectral operators,
despite the fact that the theory involves highly transcendental functions
like Airy or Bessel ones. The point is that the latter are used in the
proofs while the algorithm given at the end of Sect.~3 performs arithmetic
operations and differentiations only on {\em{explicit rational functions\/}}.

In the rest of the introduction we describe in more detail the main
results of the paper together with some of the ideas behind them.

The framework of our construction is Sato's theory of KP--hierarchy 
\cite{S, DJKM, SW, vM}. 
In particular, our eigenfunctions are Baker or wave functions
$\Psi_V(x,z)$ corresponding to planes $V$ in Sato's Grassmannian
$\Gr$ and our algebras of commuting differential operators are the spectral
algebras $\A_V$. We obtain our bispectral algebras by applying a version of
{\it
{ Darboux transformations}}, introduced in our previous paper \cite{BHY2}, on
specific wave functions which we call {\it {Bessel {\rm{(and}} Airy\/{\rm{)}}
wave functions}} (see Sect.~1 and 4). As both notions are fundamental for the
present paper we hold the attention of the reader on them.
{\it{Bessel wave functions}} are the simplest functions which solve the
bispectral problem (see \cite{Z} where they were introduced and \cite{F}). They
can be defined as follows. For $\beta = (\beta_1,\ldots ,\beta_N) \in \Cset^N$
$\Psi_\beta(x,z)$ is the unique wave function satisfying
$$ x \partial_x \Psi_\beta(x,z) = z \partial_z \Psi_\beta(x,z)$$
(i.e. $\Psi_\beta(x,z)$ depends only on $x z$) and
$$L_\beta(x,\partial_x) \Psi_\beta(x,z) = z^N \Psi_\beta(x,z),$$
where
$L_\beta(x,\partial_x) =
x^{-N}(x\partial_x -\beta_1) \cdots (x\partial_x -\beta_N)$ is the {\em{Bessel
operator}}. Obviously, the above equations lead to
$$L_\beta(z,\partial_z) \Psi_\beta(x,z) = x^N \Psi_\beta(x,z).$$
Similarly, for $\a = (\a_0, \a_1, \ldots, \a_{N-1}) \in \Cset^{N-1}$ consider
the (generalized) {\em{Airy wave function}} (see \cite{KS, Dij}) satisfying:
$$\left(\partial_x^N + \sum_{i=2}^{N-1} \alpha_i \partial_x^{N-i} - \alpha_0
x \right) \Psi_\a(x,z) = z^N \Psi_\a(x,z).$$
It depends only on $\a_0x + z^N$ and again gives a simple solution to the
bispectral problem. The Airy case is in  many respects similar to the Bessel
one. As we find the latter case richer in properties, we pay more attention to
it, contenting ourselves only with a sketch of the former.

Classically, a Darboux transformation \cite{BC, Da} of a differential 
operator $L$,
presented as a product $L = Q P$, is defined by exchanging the places of the
factors, i.e. $ \overline L = P Q $. Obviously, if $\Psi(x,\lambda)$ is an
eigenfunction of $L$, i.e. $L(x,\partial_x) \Psi(x,\lambda) = \lambda
\Psi(x,\lambda)$, then $P\Psi(x,\lambda)$ is an eigenfunction of $\overline L$.

Here we introduce Darboux transformations not only on
individual operators but also on the entire spectral algebra corresponding
to a Bessel (or Airy) plane.
In other words, we
apply them  on
operators $L$ which are polynomials $h(L_\beta)$ of Bessel (or Airy)
operators.
These transformations may be considered as B\"acklund--Darboux
transformations on the corresponding wave functions \cite{AvM}.
Such Darboux transformation is completely determined by a choice of a
$\Zset_N$-invariant operator $P(x,\partial_x)$ with rational coefficients
normalized appropriately by a factor $g^{-1}(z)$ to ensure that
$$
\Psi_W(x,z) = \frac{1}{g(z)} P(x,\partial_x) \Psi_\beta(x,z)
$$
is a wave function.
We call $ \Psi_W(x,z) $ (respectively $W$) a {\em{ polynomial Darboux
transformation}} of $\Psi_{\beta}$ (respectively $V_\beta$).
The definition of {\em{polynomial Darboux transformations\/}} of Airy planes is
similar to that in the Bessel case with only minor
modifications: $P$ is not necessarily $\Zset_N$-invariant and $g(z)$ has to
belong to $\Cset[z^N]$.

Thus we come to our main result.
\bth{0.1}
If the wave function $\Psi_W(x,z)$ is a polynomial Darboux
transformation of a Bessel or Airy wave function $\Psi_\beta(x,z)$, 
then it is a
solution to the bispectral problem, i.e. there exist differential operators
$L(x,\partial_x), \; \Lambda(z,\partial_z) $ and functions $f(z), \; \Theta(x)$
such that {\rm\eqref{0.1}} and {\rm\eqref{0.2}} hold.
\eth
Note the difference between the classical definition and the definition
introduced here. In contrast to \cite{DG} where the authors make a finite
number of ``rational'' Darboux transformations, we perform only {\em{one\/}}
polynomial Darboux transformation to achieve the same result.

Our definition of polynomial Darboux transformation is
constructive as
$P(x,\partial_x)$ is determined by the finite dimensional space 
$\ker P$.
For this reason one can explicitly present at least one operator
$L\in \A_W$; it can be given by $P h(L_\beta) P^{-1}$. Usually it
is of high order. But as it is only one element of the whole bispectral
algebra $\A_W$ there can be eventually operators of a lower order. For example,
the bispectral operators of \cite{DG} are of order two.
There is a simple procedure (see \cite{BHY2}) to produce the entire
bispectral algebra $\A_W$ of commuting differential operators.
In addition, one can show that the {\em{spectral curve\/}} $ \spec \A_W$ (see
e.g. \cite{AMcD} for definition) {\em{is rational, unicursal and
$\Zset_N$--invariant}}.

In the course of our work we have widely used important ideas introduced by
G.~Wilson \cite{W}. Among them we mention first the idea of explicitly writing
conditions on vectors of a plane $V \in \Gr$ which define the new plane
obtained by a Darboux transformation. Second is the notion of 
involutions on the Sato's Grassmannian. In particular, we extend the
{\em{bispectral involution\/}} $b$ introduced in \cite{W} to the manifolds of
polynomial Darboux transformations. More precisely, we prove the following 
theorem, from which \thref{0.1} is an obvious consequence.

\bth{0.2}
{\em{(i)}} The bispectral involution
is defined for planes $W$ which are polynomial Darboux transformations
of Bessel or Airy planes.

{\em{(ii)}} The image $b W$ of such a plane $W$ is again a polynomial Darboux
transformation
of the corresponding Bessel (respectively) Airy plane.
\eth

Our main concern in the present paper is to prove \thref{0.2}. Our second goal 
is to provide explicit formulae and examples (see Sect. 5),
which are not only an illustration of our method but also show the existence
of new families of bispectral operators with particular properties. Some of
them generalize directly the well known ones like the Duistermaat--Gr\"unbaum's
{\em{``even case''\/}} \eqref{0.6} \cite{DG}. Other families exhibit quite
different properties from
the well known examples. In this respect Sect.~5 has also the role to supply
diverse experimental material for new insights into the theory of bispectral
algebras. We draw the attention of the reader also to the explicit
formulae for the action of the bispectral involution on an important class of
Darboux transformations (which we call monomial) of Bessel operators.
As a particular case, we obtain such formulae for all second order bispectral 
operators found in \cite{DG}. 

The class of monomial Darboux transformations has also other remarkable
properties, e.g. they are connected to representation theory of
$W_{1+\infty}$--algebra. We do not touch this matter here for lack of space. 
The interested reader can learn about it in \cite{BHY3}.

A natural question is if the operators found in this paper form the entire
class of bispectral operators. The answer is negative as recently shown
in \cite{BHY4}.

At the end for the reader's convenience we give a brief description
of the organization of the paper. Sect.~1 is intended only for reference.
It reviews results connected with Sato's theory, which we need for the 
treatment
of the bispectral problem. Besides the general notions (see e.g. 
\cite{S, DJKM, SW}) we recall the involutions, introduced by G.~Wilson
\cite{W} and in particular, the bispectral involution.
In Sect.~2 we introduce our manifolds $\GRBN$ of polynomial
Darboux transformations of Bessel planes. We give two equivalent definitions
(\deref{5.5} and the one provided by the statement of \thref{5.7}). Sect.~3
contains our main results -- 
Theorems \ref{t0.1} and \ref{t0.2}. for the Bessel
case. Sect.~4 deals with the analogs of Sect.~2 and 3 for the Airy case 
(although
in different order). The last Sect.~5 is devoted to explicit examples of 
bispectral
operators, which have been studied in other papers \cite{DG, W2},
as well as new families (which we have not seen elsewhere). The emphasis in 
Sect.~5
is rather on the simple algorithmic way of constructing bispectral operators
(wave functions, etc.) than on the novelty of the examples.

For readers who wish to see the main results as soon as possible we propose
another plan of reading the paper. They can start with Sect.~2 and read 
it up to the statement of \thref{5.7}, returning to Sect.~1 for reference
when needed. Then skipping the (technical) proof of \thref{5.7}, they can go
Sect.~3. After that, taking for granted the proof of \thref{6.2},
they can look at the examples 
of bispectral operators, originating from Bessel ones in Sect.~5. Thus they
will have a complete picture of the results in the Bessel case,
and having this experience, they can easily go through the Airy case.

More detailed information about the material included in each section
can be found in its beginning.

The present paper is a part of our project on the bispectral problem
\cite{BHY1}--\cite{BHY4}. The main results contained here were announced at the
conference of Geometry and Mathematical Physics, Zlatograd 95 
(see \cite{BHYc}).
\hfill\\

After this paper was written, we got a paper \cite{KRo} 
where some of the results
about the Airy case were obtained independently.

{\flushleft{\bf{Acknowledgements}}}

\medskip\noindent

We are grateful to F.~A.~Gr\"unbaum and G.~Wilson for their interest 
in the paper and for suggestions which led to improving the presentation
of our results. We also thank the referee who proposed important changes
towards making the text more
``reader friendly''.
This work was partially supported by Grant MM--523/95 of Bulgarian
Ministry of Education, Science and Technologies.

\sectionnew{Preliminaries}
In this section we have collected results about Sato's theory, relevant to 
the bispectral problem. For reader's convenience we have divided the section
into 4 subsections, whose titles, hopefully, give an idea of their content.
The reader, who is acquainted with Sato's theory may even skip this section and
return to it for references when needed. More detailed account of the 
material of the subsections can be found in their beginnings.
\subsection{Sato's Grassmannian and KP--hierarchy}
We shall recall some facts and notation from Sato's theory
of KP-hierarchy needed in the paper. The survey below cannot be used as
a systematic study. There are several complete texts on Sato's theory, starting
with the original papers of M. Sato and his collaborators \cite{S, DJKM} (see
also \cite{SW, vM}).

Consider the space of formal series
$$
\Vset=\Bigl\{ \sum_{k\in\Zset} a_kv_k \Big|\; a_k=0\ {\rm for}\
k\gg 0\Bigr\}.
$$
{\em Sato's Grassmannian\/} $Gr$ \cite{S, DJKM, SW} consists of all
subspaces ({\em{planes\/}}) $W\subset\Vset$ which have an admissible basis
$$
w_k=v_k+\sum_{i<k}w_{ik}v_i,\quad k=0,1,2,\ldots
$$

In Sato's theory $\Vset$ is most often realized as the space of formal Laurent
series in $z^{-1}$ via $v_k=z^k$.
The {\em{Baker\/}} (or {\em{wave function\/}}) $\Psi_W(x,z)$ of the plane $W$
contains the whole information about $W$  as the vectors
$w_{k}=\partial^k_x\Psi_W(x,z)|_{x=0}$
form an admissible basis of $W$.
We can expand $\Psi_W(x,z)$ in a formal series
\beq
\Psi_W(x,z)=e^{x z}
\left(1+\sum_{k=1}^\infty a_k(x)z^{-k}\right).
\label{1.13}
\eeq
The wave function $\Psi_W(x,z)$ can also be written in terms of the so-called
wave operator $K_W$. This is a pseudo-differential operator defined by
\beq
K_W(x,\partial_x) =1+\sum_{k=1}^\infty a_k(x)\partial_x^{-k}.
\label{1.14}
\eeq
Then obviously
$$
\Psi_W(x,z)=K_W(x,\partial_x) e^{x z}.
$$
Introduce also the pseudo-differential operator
\beq
P=K_W \circ\partial_x\circ K_W^{-1}.
\label{1.15}
\eeq
For the treatment of the bispectral problem the following identity is crucial:
\beq
P\Psi_W(x,z)=z\Psi_W(x,z).
\label{1.17}
\eeq
When it happens that some power of $P$, say $P^N$, is a differential operator
$L$ we get that $\Psi_W(x,z)$ is an eigenfunction of an ordinary differential
operator:
\beq
L\Psi_W(x,z) = z^N\Psi_W(x,z).
\label{1.18}
\eeq
It is easy to show that $P^N$ is differential iff the plane $W$ is invariant
under the multiplication by $z^N$:
\beq
z^NW\subset W.
\label{1.19}
\eeq
The submanifold of $Gr$ consisting of planes $W$ satisfying
\eqref{1.19} is denoted by $Gr^{(N)}$.

A very important  object connected to the plane $W$ is the
algebra $A_W$ of polynomials $f(z)$ that leave $W$ invariant:
\beq
A_W=\{ f(z) | f(z)W\subset W\}.
\label{1.20}
\eeq
For each $f(z)\in A_W$ one can show that there exists a unique
differential operator $L_f(x, \partial_x)$, the order of $L_f$ being equal to
the degree of $f$, such that
\beq
L_f(x, \partial_x)\Psi_W(x,z)=f(z)\Psi_W(x,z).
\label{1.20'5}
\eeq
Explicitly we have
\beq
L_f=K_W\circ f(\partial_x)\circ K_W^{-1}.
\label{1.21}
\eeq
We denote the commutative algebra of these operators by $\A_W$, i.e.
\beq
\A_{W} = \{ L_{f} | L_{f} \Psi_{W} = f \Psi_{W}, \; f\in A_{W} \}.
\label{1.22}
\eeq
Obviously, $A_W$ and $\A_W$ are isomorphic.
We call $\A_W$ {\em spectral algebra\/} corresponding to the plane $W$.
Following I.~Krichever (see e.g.
\cite{KrN}), we introduce a {\em  rank\/} of $\A_W$ to be the dimension of the
space of eigenfunctions $\Psi_W$; this number is equal to
the greatest common divisor of the orders of the operators $L_f$.
At the end we define the {\em spectral curve\/}
corresponding to the plane $W$ to be
$\spec\A_W$ (for definition see e.g. \cite{AMcD}).
It is known that $\spec \A_W$ is an algebraic curve
(see \cite{BC, I, KrN}).
\bre{1.2}
If $\Psi_W(x,z)$ is well defined for $x=x_0$ we set $v_k=e^{x_0 z} z^k$
and consider the subspace $W^{x_0}$ of $\Vset$ with an admissible basis
$w_{k}=\partial^k_x\Psi_W(x,z)|_{x=x_0}$.
The wave functions of $W^{x_0}$ and $W$ are connected by
$
\Psi_{W^{x_0}}(x,z) = e^{-x_0 z} \Psi_{W} (x + x_0, z)
$
and obviously
$$
K_{W^{x_0}}(x,\partial_x) = K_{W}(x+x_0,\partial_x).
$$
These shifts are inessential for the bispectral problem and for our proofs.
Throughout the paper we shall sometimes work with $\Psi_{W^{x_0}}$ calling it
by abuse of notation a wave function of $W$ and denoting it $\Psi_W$.
\ere
\bre{1.3}
In Sato's theory of KP hierarchy one usually considers the wave function
$$
\Psi_W(t, z)=e^{\sum_{k=1}^\infty t_k z^k}
\left(1 + \sum_{k=1}^{\infty} a_k(t) z^{-k}\right)
$$
depending on all times $t_1,t_2,\ldots$ (here $t_1 =x$). Then
the operators $K_W(t, \partial_x)$ and $P(t, \partial_x)$ are given again by
formulae (\ref{1.14}, \ref{1.15}) with $a_k(x)$ replaced with $a_k(t)$  and
$P$ satisfies the following infinite system of non-linear differential
equations
\beq
\frac{\partial}{\partial t_k} P=[P^k_+, P],
\label{1.16}
\eeq
where $P^k_+$ stands for the differential part of the $k$-th power $P^k$ of the
operator $P$, $[\,.\,,\,.\,]$ is the standard commutator  of
pseudo-differential
operators. The equations \eqref{1.16} are called the {\em KP hierarchy}. If
the plane $W$ lies in $\GRN$ then \eqref{1.16} can be written in the form
\beq
\frac{\partial}{\partial t_k} L=[L^{k/N}_+, L]
\label{1.16aa}
\eeq
($L= P^N$) and is called an {\em{$N$-th reduction\/}}
of the KP-hierarchy or an {\em{$N$-th Gelfand--Dickey hierarchy\/}} \cite{GD}.
\ere
\subsection{Darboux transformations}
In this subsection we recall the notion of {\em{Darboux transformations\/}}
on objects connected to points of Sato's Grassmannian, introduced in our recent
paper \cite{BHY2}.

Recall that a {\em{Darboux transformation\/}} \cite{Da} of an ordinary
differential operator $L$ is obtained by presenting it as a product and
exchanging the places of the factors:
$$
L = Q P \mapsto {\overline L} = P Q.
$$
A (monic) operator $L$ is completely determined by its kernel:
if $\{\Phi_0,\ldots,\Phi_{n-1}\}$ is a basis of $\ker L$ then
(see e.g. \cite{I})
\beq
L\Phi=\frac{\Wr(\Phi_0,\ldots,\Phi_{n-1},\Phi)}{\Wr(\Phi_0,\ldots,\Phi_{n-1})}
\label{3.1}
\eeq
where $\Wr$ denotes the Wronski determinant.
The next lemma answers the question when the factorization $L=QP$ is possible
(see e.g. \cite{I}).
\ble{3.1}
{\em(i)}  For a given basis $\{\Phi_0,\ldots,\Phi_{n-1}\}$ of\/
$\ker L$ set
\beqa
q_1(x)&=&\partial_x\log \Phi_0, \hfill
\label{3.2}\\
q_k(x)&=&\partial_x\log
\frac{\Wr(\Phi_0,\ldots,\Phi_{k-1})}{\Wr(\Phi_0,\ldots,\Phi_{k-2})},
\quad 2\le k\le n. \hfill
\label{3.2'5}
\eeqa
Then the operator $L$ can be factorized as follows
\beq
L=(\partial_x  -q_n)(\partial_x-q_{n-1})\cdots(\partial_x-q_1).
\label{3.3}
\eeq
{\em(ii)} $L$ can be factorized as
\beq
L=QP\ \ {\rm iff}\ \ \ker P\subset \ker L.
\label{3.4}
\eeq
In this case
\beq
\ker Q = P(\ker L).
\label{3.5}
\eeq
\ele
A slightly more general construction is the following one.
For operators $L$ and $P$ such that the kernel of $P$ is invariant under $L$,
i.e.
\beq
L(\ker P)\subset \ker P
\label{3.6}
\eeq
we consider the transformation
\beq
L\mapsto \overline L=PLP^{-1}.
\label{3.7}
\eeq
The fact that $\overline L$ is a differential operator follows from \leref{3.1}
(ii). Indeed, $L(\ker P)\subset \ker P$ is equivalent to $\ker P\subset
\ker(PL)$.

In \cite{BHY2} we defined a version of Darboux transformation on points in
Sato's Grassmannian and on related objects -- wave functions, tau-functions and
spectral algebras.
\bde{3.2}
We say that a plane $W$ (or the corresponding wave
function $\Psi_W(x,z)$) is a {\em Darboux transformation\/} of the plane $V$
(respectively wave function $\Psi_V(x,z)$) iff there exist monic polynomials
$f(z)$, $g(z)$ and differential operators $P(x,\partial_x)$, $Q(x,\partial_x)$
such that
\beqa
&&\Psi_W(x,z)=\frac{1}{g(z)} P(x,\partial_x) \Psi_V(x,z),
\label{3.8} \\
&&\Psi_V(x,z)=\frac{1}{f(z)} Q(x,\partial_x) \Psi_W(x,z).
\label{3.9}
\eeqa
\ede
An equivalent definition is that $W$ is a Darboux transformation of $V$ iff
\beq
fV\subset W\subset \frac1gV
\label{3.13}
\eeq
for some polynomials $f(z)$, $g(z)$.

Simple consequences of \deref{3.2} are the identities
\beqa
&&PQ\Psi_W(x,z)=f(z)g(z)\Psi_W(x,z),
\label{3.10}\\
&&QP\Psi_V(x,z)=f(z)g(z)\Psi_V(x,z).
\label{3.11}
\eeqa
The operator $\overline L=PQ\in\A_W$ is a Darboux transformation of
$L=QP\in\A_V$.

Having in mind applications to the bispectral problem, the most important for
our study is the case when
\beq
A_V=\Cset[z^N],\quad \A_V=\Cset[L_V]
\label{3.14}
\eeq
for some natural number $N$ and a differential operator $L_V$ of order $N$.
(This is the simplest case of rank $N$ spectral algebra with a spectral curve
$\Cset$.) Then due to \eqref{3.11} we have
\beqa
&&f(z)g(z)=h(z^N),
\label{3.15}\\
&&QP=h(L_V)
\label{3.16}
\eeqa
for some polynomial $h(z)$.
In \cite{BHY2} we connected the spectral algebra $A_W$ (respectively $\A_W$)
with $A_V$ (respectively $\A_V$).
\bpr{3.5}
{\rm{(i)}}
The Darboux transformations preserve the rank of the spectral algebras, i.e.\
if $W$ is a Darboux transformation of $V$ then $\rank \A_W = \rank \A_V$.

{\rm{(ii)}}
If $\A_V=\Cset[L_V]$, $\ord L_V = N$ then
\beq
A_W=\left\{ u\in \Cset[z^N] \mid u(L_V)\ker P\subset \ker P\right\}
\label{3.18}
\eeq
and
\beq
\A_W=\left\{ Pu(L_V)P^{-1} \mid u\in A_W\right\}.
\label{3.19}
\eeq
\epr
\subsection{Bessel operators, Bessel planes and related objects}
Now we define the planes of the Sato's Grassmannian on which we shall perform
the Darboux transformations.
 For $\beta \in \Cset^N$ such that
\beq
\sum_{i=1}^{N}\beta_i = \frac{N(N-1)}{2}
\label{2.15}
\eeq
we introduce the ordinary differential operator
\beq
  P_{\beta}(D_z) = (D_z-\beta_1)(D_z-\beta_2) \cdots (D_z-\beta_N),
  \label{2.13}
\eeq
where $D_z = z\partial_z$, and consider the differential equation
\beq
 P_{\beta} (D_z) \Phi_{\beta}(z) = z^N \Phi_{\beta} (z).
 \label{2.14}
\eeq
For every sector $S$ with a center at the irregular singular point $z= \infty$
and an angle less than $2 \pi$ the equation \eqref{2.14} has a solution
$\Phi_{\beta}$ with an asymptotics
\beq
 \Phi_{\beta}(z) \sim\Psi_{\beta}(z) =
e^z \Bigl( 1+ \sum_{k = 1}^{\infty} a_k (\beta) z^{-k} \Bigr)
\label{2.14'5}
\eeq
for $|z| \to \infty$, $z \in S$ (see e.g. \cite{Wa}). Here
$a_k(\beta)$ are symmetric polynomials in $\beta_i$.
The function $\Phi_\beta(z)$ can be taken to be (up to a rescaling) the Meijer's
$G$-function
\beq \label{1.33'2}
\Phi_\beta(z) = G^{N0}_{0N}\Bigl(\Bigl(\frac{-z}{N}\Bigr)^N \Big|
\frac1{N}\beta\Bigr)
\eeq
-- see \cite{BE}, \S $5.3$.

The next definition is fundamental for the present paper.
\bde{2.2}
{\it {Bessel wave function}} is called the function
$ \Psi_{\beta}(x,z) = \Psi_{\beta}(xz) $ (cf. \cite{F, Z}).
The {\it {Bessel operator}} $L_{\beta}$ is defined as
\beq
 L_{\beta}(x,\partial_x) = x^{-N} P_{\beta} (D_x).
 \label{2.16}
\eeq
A Bessel wave function $\Psi_\beta$ defines a plane $V_\beta \in\Gr$ (called
{\em{Bessel plane\/}}) by the standard procedure:
$$ V_\beta = \span \{ \partial_x^k \Psi_\beta(x,z)|_{x = 1} \} .$$
(In fact $\Psi_{V_\beta}(x,z) = e^{-z} \Psi_\beta(x+1,z)$, cf.\
\reref{1.2}; we took $x_0=1$ because $\Psi_\beta(x,z)$ is singular at $x=0$,
arbitrary $x_0\not=0$ will do.)
\ede
Because $\Psi_\beta(x,z)$ depends only on $xz$, i.e.
\beq
D_x \Psi_\beta(x,z)= D_z \Psi_\beta(x,z),
\label{2.17'1}
\eeq
it gives the simplest solution to the bispectral problem:
\beqa
&&L_\beta (x, \partial_x) \Psi_\beta(x,z) =
   z^N \Psi_\beta(x,z),
\label{2.17'2}\\
&&L_\beta (z, \partial_z) \Psi_\beta(x,z) =
   x^N \Psi_\beta(x,z).
\label{2.17'3}
\eeqa
\subsection{Involutions in Sato's Grassmannian}
In this subsection, following G.~Wilson \cite{W}, we define several involutions
on points of Sato's Grassmannian and  on related objects -- wave functions and
wave operators. Besides the general properties of the
involutions taken from \cite{W}, we specify their action on Bessel planes.

Introduce after \cite{DJKM} the non-degenerate form in $\Vset$ (realized as the
space of formal Laurent series in $z^{-1}$)
$$
B(f,g)=-\res_\infty f(z)g(-z)\,dz, \quad f,g\in\Vset.
$$
If $V\in Gr$ is a plane, define $aV\in  Gr$ to be the plane orthogonal to $V$
with respect to the form $B(\,.\,,\,.\,)$, to wit
\beq
aV=\{ g(z) \mid B(f,g)=0,\ {\rm for\ all}\ f\in V\}.
\label{4.1}
\eeq
Obviously, $a(aV)=V$, i.e.\ the map $a$ is an involution. Following \cite{W},
we call it the {\it  adjoint involution}.
On the wave operator $K_W$ \eqref{1.14} the involution $a$ acts as \cite{W}
\beq
K_{aV}=\left(K^*_V\right)^{-1},
\label{4.3}
\eeq
 where $*$ is the formal conjugation on pseudo-differential operators, i.e.\
the antiautomorphism defined by $\partial^*_x=-\partial_x$, $x^*=x$.
For our purposes the most important property of the involution $a$ is that it
inverses inclusions, i.e.:
\beq
{\rm if}\ W\subset V,\ {\rm then}\ aW\supset aV.
\label{4.4}
\eeq
The following proposition will be used in the description of the action of $a$
on Darboux transformations. Its simple proof is similar to that of Corollary
7.7 from \cite{W}.
\bpr{4.3}
{\rm{(i)}} If
$
\Psi_W(x,z)=\dfrac{1}{g(z)} P(x,\partial_x) \Psi_V(x,z),
$
then
\beq
\Psi_{aV}(x,z)=\frac{1}{\check g(z)} P^*(x,\partial_x) \Psi_{aW}(x,z)
\label{4.7}
\eeq
where $\check g(z)=g(-z)$.

{\rm{(ii)}} {\rm{\cite{W}}} Let $\A_V$ be the algebra of operators
\eqref{1.22}.  Then $\A_{aV}$ consists of the conjugated operators of $\A_V$.
\epr
\noindent
\proof
We have
$
\Psi_W(x,z) = g^{-1}(z) P\Psi_V(x,z) = g^{-1}(z) P K_Ve^{xz}=
P K_V g^{-1}(\partial_x)e^{xz}
$
which implies
$K_W=P K_V g^{-1}(\partial_x).$
Applying the involution $a$ we obtain
$$
K_{aW}=\left(K^*_W\right)^{-1}=(P^*)^{-1} \left(K^*_V\right)^{-1}
\check g(\partial_x),
$$
yielding
$ P^*K_{aW}=K_{aV} \check g(\partial_x) $
and hence \eqref{4.7}.
\qed

\smallskip

The {\em  sign involution\/} $s$ \cite{W} is defined on the wave functions
by the property
\beq
\Psi_{sV}(x, z)=\Psi_V(-x,-z).
\label{4.8}
\eeq
On the wave operators $K_V(x,\partial_x)$ \eqref{4.8} translates into
\beq
K_{sV}(x,\partial_x)=K_V(-x,-\partial_x).
\label{4.9}
\eeq
The plane $sV$ is defined by
\beq
sV=\{f(-z) \mid f(z)\in V\}.
\label{4.10}
\eeq

In the subsequent chapters we shall need the action of $a$ and $s$ on the
Bessel planes. We describe them in the next proposition.
\bpr{4.4}
 The involutions $s$ and $a$ act on Bessel planes
$V_\beta$ {\rm(}$\beta \in \Cset^N${\rm)} as follows
\beqa
&&sV_\beta=V_\beta,
\label{4.12}\\
&&aV_\beta=V_{a(\beta)},
\label{4.13}
\eeqa
 where $a(\beta)=(N-1)\delta-\beta$, $\delta = (1,1,\ldots,1).$
\epr
\noindent
\proof
We compute the action of the involutions on the corresponding
wave functions. Obviously,
$$
\Psi_{sV_\beta}(x,z)=\Psi_{\beta}(-x,-z)=\Psi_{\beta}(x,z),
$$
showing \eqref{4.12}.

\prref{4.3}, eq.\ \eqref{2.17'2} and the fact that
$
L^*_\beta=(-1)^N L_{a(\beta)}
$
imply
$
\Psi_{a V_\beta}(x,z) = \gamma(z) \Psi_{V_{a(\beta)}}(x,z)
$
for some formal power series $\gamma(z)$ in $z$. To show that $\gamma(z)=1$ we
notice that $\Psi_{a V_\beta}(x,z)$ depends on $xz$. Indeed, \eqref{2.17'1}
implies
$
K_{V_\beta}(cx, c^{-1}\partial_x) = K_{V_\beta}(x, \partial_x)
$
for all $c \not = 0$ and the same is true for
$K_{a V_\beta} = (K_{V_\beta}^*)^{-1}.$
\qed

\smallskip

We end this section by recalling the {\it bispectral involution\/} $b$ which
Wilson \cite{W} introduced for the purpose of the bispectral problem.
Contrary to $a$ and $s$, 
the bispectral involution $b$ {\it is not defined\/} on
the entire Grassmannian. Whenever one can define $b$, put
\beq \label{1.42'5}
\Psi_{bV}(x,z)=\Psi_V(z,x),
\eeq
i.e.\ the involution $b$ interchanges the places of the arguments $x$ and $z$.

A simple example of a point $V\in Gr$ where the involution $b$ is well defined
is the Bessel plane $V_\beta$. It immediately follows 
 {} from the definition that
$$
\Psi_\beta(x,z)=\Psi_\beta(z,x),
\quad
{\rm i.e.}\ bV_\beta=V_\beta.
$$

In terms of the bispectral involution our approach to the bispectral problem
can be formulated geometrically as follows:

\medskip

Find points $V\in Gr$ such that

1) $g(z)V\subset V$  for some nontrivial polynomial $g(z)$;

2) $bV$ exists and $f(z)bV\subset bV$ for some nontrivial polynomial $f(z)$.

\medskip
A very important general property of $b$, which we intend to use, is its
connection to the other involutions \cite{W}:
\beq
ab=bas.
\label{4.18}
\eeq
For completeness we also point out that the involution $s$ commutes with $a$
and $b$:
\beq
as=sa,\quad bs=sb.
\label{4.19}
\eeq
\sectionnew{Polynomial Darboux transformations of\\ Bessel wave functions}
The main purpose of this section is to introduce the submanifolds
(denoted below by $\GRBN$) on which, as we prove in the next section,
the bispectral involution $b$ 
is well-defined, and whose points correspond to bispectral operators.
The points of $\GRBN$ are obtained by a version of B\"acklund--Darboux
transformation performed on Bessel wave functions (or equivalently,
on polynomials $h(L_\b)$ of Bessel operators). Below we call these
transformations {\em polynomial Darboux transformations\/}.
\deref{5.5}, where this is done and the statement of \thref{5.7},
where we provide an equivalent definition, form the heart of 
the present section. \deref{5.5} has the advantage to be more natural 
and to supply an algorithmic procedure for constructing bispectral
operators. The second definition (from \thref{5.7}) is more suitable
for the proof of our bispectrality theorem in Sect.~3.

The reader who wishes to see as soon as possible the main results
of the paper can use the second definition, the one from the
statement of \thref{5.7}, skipping its technical proof, which
occupies half of the section.

In the first half of this section we describe the kernel of the
operator $P$ from \deref{3.2} and (which is equivalent) -- the 
conditions of the type as in \cite{W}, imposed on a Bessel plane,
which define the corresponding Darboux transformation.
To do so, we first need a description of the kernels of the operators
$h(L_\b)$ which are polynomials of $L_\b$.

Fix $\beta\in\Cset^N$ satisfying \eqref{2.15} and let $V_\beta$ be the
corresponding Bessel plane (see Subsect.~1.3). Throughout this section $W$ will
be a Darboux transformation of $ V_\beta$ (we shall use the notation of
\deref{3.2} with $V=V_\beta$). We shall need a lemma describing the kernel of
the operator $h(L_\beta)$ for an arbitrary polynomial $h$.
\ble{5.1}
 Let $h(z)$ be a polynomial
\beq
h(z)=z^{d_0}\left(z-\lambda_1^N\right)^{d_1}\cdots
\left(z-\lambda_r^N\right)^{d_r},\quad
\lambda_i^N\ne \lambda_j^N,\ \lambda_0=0,\ d_i\ge0.
\label{5.1}
\eeq
Then we have

{\em(i)} $\ker h(L_\beta)=\bigoplus_{i=0}^r \ker
\left(L_\beta-\lambda_r^N\right)^{d_i}$.

{\em(ii)} $(L_\beta)^d=L_{\beta^d}$,  where
\beq
\beta^d=(\beta_1,\beta_1+N,\ldots,\beta_1+(d-1)N,\ldots,\beta_N,\ldots,\beta_N
+(d-1)N).
\label{5.1'5}
\eeq

{\em(iii)}  If
$\{\beta_1,\ldots,\beta_N\}=\{\underbrace{\alpha_1,\ldots,\alpha_1}_{k_1}
,\ldots,
\underbrace{\alpha_s,\ldots,\alpha_s}_{k_s} \}$
with distinct $\alpha_1, \ldots,\alpha_s$, then
$$
\ker L_\beta=\span\left\{ x^{\alpha_i}(\ln x)^k\right\}_{1\le i\le s,\ 0\le
k\le k_i-1}.
$$

{\em(iv)}  For $\lambda\ne 0$
$$
\ker\left(L_\beta-\lambda^N\right)^d =\span\left\{
\partial_z^k\Psi_\beta(x,z)|_{z=\lambda\varepsilon^j}
\right\}_{0\le k\le d-1,\
0\le j\le N-1},
$$
where $\varepsilon=e^{2\pi i/N}$ is an $N$-th root of unity.
\ele
The {\em proof\/} being obvious is omitted (cf.\ \leref{3.1}).
\qed

Let us consider the simplest {\bf {example}} of a Darboux transformation.
Set
\beq
h(z)=z^d,\quad g(z)=z^n,\quad f(z)=z^{dN-n}
\label{5.2}
\eeq
and
$\gamma=\beta^d$, i.e.\
\beq
\gamma_{(k-1)d+j} := \beta_k+(j-1)N,\quad 1\le k\le N,\ 1\le j\le d.
\label{5.3}
\eeq
For an $n$-element subset $I$ of $\{1,\ldots, dN\}$
such that $\gamma_i\not=\gamma_j$ for $i\not=j\in I$, we put
\beq
\ker P=\span\bigl\{x^{\gamma_i}\bigr\}_{i\in I}.
\label{5.5}
\eeq
Such $P$ corresponds
to a Darboux transformation $\Psi_I(x,z)$ of $\Psi_\beta(x,z)$.
The following simple fact will be useful in the sequel.
\ble{5.2}
$\Psi_I(x,z)$  is again a Bessel wave function{\rm:}
\beq
\Psi_I(x,z)=\Psi_{\gamma+dN\delta_I -n\delta} (x,z).
\label{5.6}
\eeq
Here and further we use the vectors $\delta_I$, $\delta$ defined by
\beq
(\delta_I)_i=\cases{1, &if $i\in I$\cr
                    0, &if $i\not\in I$\cr}
\label{5.7}
\eeq
 and
\beq
\delta_i=1\ \ {\it for\ all}\ \ i\in\{1,\ldots, dN\}.
\label{5.8}
\eeq
\ele
\noindent
{\it Proof}. By definition
\beq
\Psi_I(x,z)=z^{-n} L_{\gamma_I}\Psi_\gamma(x,z),
\label{5.9}
\eeq
where $\gamma_I=\{\gamma_i\}_{i\in I}$. Then $\Psi_I(x,z)$ is an eigenfunction
of the differential operator
$L_{\gamma_I} L_\gamma L_{\gamma_I}^{-1}$, which is
straightforwardly computed to be equal to $L_{\gamma+dN\delta_I-n\delta}$.
\qed

\smallskip

Our next step is to study the spectral algebra $\A_\beta\equiv \A_{V_\beta}$ of
a Bessel plane $V_\beta$ (see (\ref{1.20}, \ref{1.22})).
\ble{5.3}
 If $L(x,\partial_x)\Psi_\beta(x,z)=u(z)\Psi_\beta(x,z)$
for some operator $L\in\A_\beta$ and some polynomial $u(z)\in A_\beta$, then
$L$ is a linear combination of Bessel operators $L_\alpha$,
$\alpha\in\Cset^k$ such that $L_\alpha\Psi_\beta(x,z)= z^k\Psi_\beta(x,z)$.
\ele
\noindent
{\it Proof}. Let $u(z)=\sum u_kz^k$, $u_k\not=0$, $0\le k\le M$ for some $M$.
Then for arbitrary $c\ne0$  we have
$$
u(cz)\Psi(x,z)=u(cz)\Psi(c^{-1}x,cz)=
L(c^{-1}x, c\partial_x)\Psi_\beta(x,z).
$$
This implies that $u(cz)\in A_\beta$
and thus $z^k\in A_\beta$. On the other hand $D_zV_\beta\subset V_\beta$ and
the compatibility condition is of the form $(D_z-\alpha_1)\cdots
(D_z-\alpha_k)\Psi_\beta(1,z)=z^k\Psi_\beta(1,z)$ which implies $L_\alpha
\Psi_\beta(x,z)=z^k\Psi_\beta(x,z)$.
\qed

\smallskip

Let us introduce the following terminology. We say that $\beta\in\Cset^N$ is
{\it generic\/} if $V_\beta$ is not a Darboux transformation of another Bessel
plane $V_{\beta'}$ with $\beta'\in \Cset^{N'}$, $N'<N$.
The next proposition seems obvious but we do not know a simpler proof.
\bpr{5.4}
 For a generic $\beta\in \Cset^N$ we have
\beq
A_\beta=\Cset[z^N], \qquad \A_\beta = \Cset[L_\beta].
\label{5.10}
\eeq
\epr
\noindent
{\it Proof}. We shall prove that if $\rank A_\beta=r<N$ then there exists
$\beta'\in
\Cset^r$ (with $A_{\beta'}=\Cset[z^r]$) such that $V_\beta$ is a Darboux
transformation of $V_{\beta'}$. The main idea is to apply to $V_\beta$ Darboux
transformations which lead again to Bessel planes but reduce the order of the
operator $L_\beta$. Note that according to \prref{3.5} they preserve the
rank $r$ of the spectral algebra.

Split the set $\beta$ into congruent $\mod N\Zset$ classes
$$
\beta=\left(\beta^{(1)}_1,\ldots, \beta^{(1)}_{N_1},\ldots , \beta_1^{(p)},
\ldots, \beta_{N_p}^{(p)}\right)
$$
such that
$$
\beta_s^{(i)}-\beta_t^{(j)}\not\in N\Zset \quad{\rm for}\ \ i\ne j\ \ {\rm
and\ all}\ \ s,t
$$
($N=N_1+\cdots+ N_p$).

By a Darboux transformation this $\beta$ can be changed to
$$
\beta'=\bigl( \underbrace{\beta^{(1)},\ldots,\beta^{(1)}}_{N_1}
,\ldots,
\underbrace{\beta^{(p)},\ldots,\beta^{(p)}}_{N_p}
\bigr)
$$
such that
\beq
|\Re(\beta^{(i)}-\beta^{(j)})|<N \quad {\rm and} \ \ \beta^{(i)}\ne \beta^{(j)}
\ \ {\rm for}\ \ i\ne j
\label{5.11}
\eeq
(see \leref{5.2}).

Suppose that $A_{\beta'}\ne \Cset[z^N]$. Then by \leref{5.3} there exists a
Bessel operator $L_\alpha$ such that
\beq
L_\alpha\Psi_{\beta'}(x,z)=z^M\Psi_{\beta'}(x,z),\quad \alpha\in \Cset^M
\label{5.12}
\eeq
and $L_\alpha\ne L^k_{\beta'}$ for any $k$. It is clear that
\beq
L_\alpha L_{\beta'}=L_{\beta'}L_\alpha,
\label{5.13}
\eeq
which is equivalent to
$$
\left\{ \alpha_1+N,\alpha_2+N,\ldots, \alpha_M+N, \beta^{(1)},\ldots,
\beta^{(1)}, \ldots,\beta^{(p)},\ldots,\beta^{(p)}\right\}
$$
$$
{}=\left\{ \beta^{(1)}+M,\ldots, \beta^{(1)}+M, \ldots, \beta^{(p)}+M, \ldots,
\beta^{(p)}+M, \alpha_1,\ldots,\alpha_M\right\}.
$$
Now if $M>N$ this imply that $\beta'\subset\alpha$ and
therefore there exists a Bessel operator $L_{\alpha'}$ such that
$$
L_\alpha=L_{\alpha'}L_{\beta'}
\quad{\rm and}\quad
L_{\alpha'}L_{\beta'}=L_{\beta'}L_{\alpha'}.
$$
Repeating the same argument with $\alpha'$, we obtain that there exists
$L_\alpha $ satisfying \eqref{5.13} with $M<N$. But then \eqref{5.13} is
equivalent to $V_{\beta'}=V_\alpha$. By \prref{3.5}
$r=\rank A_\beta=\rank A_{\beta'}=\rank A_\alpha$ divides $M$ and $N$.
If $V_\alpha=V_{\beta'}= \Cset [z^r]$ this finishes the proof. Otherwise we can
repeat the above argument with  $V_{\alpha}$ instead of $V_{\beta'}$.
\qed

\smallskip

Now we come to the main purpose of this section: the definition of manifolds of
Darboux transformations, which will give solutions to the bispectral problem.
To get some insight we shall consider, following Wilson \cite{W}, the
geometrical meaning of Darboux transformations, provided by the so-called
conditions $C$.

\prref{5.4} implies that for generic $\beta\in\Cset^N$ (\ref{3.15},
\ref{3.16}) hold with $V=V_\beta$
and $\ker P$ is a subspace of $\ker h(L_\beta)$. Each element $f$ of
$\ker P$ corresponds to a {\em{condition\/}} $c$ (a linear functional
on $V_\beta$), such that
\beq
f(x)=\langle c,\Psi_\beta(x,z)\rangle,
\label{5.16}
\eeq
$c$ acts on the variable $z$.
These linear functionals form an $n$-dimensional linear space $C$
({\em{space of conditions\/}}) where $n=\ord P$. In this terminology
the definition of
Darboux transformation can be reformulated as
$$
W = \frac{1}{g(z)} \Bigl\{ v\in V_\beta\; \Big| \;\langle c,v \rangle
= 0  \;
{\textrm{ for all }} c\in C \Bigr\}
$$
(see  \cite{W, BHY2}).
Following Wilson \cite{W}, we call the condition $c$ supported at $\lambda$
iff it is of the form (cf.\ \leref{5.1} (iv))
\beq
c=\sum_k a_k\partial_z^k|_{z=\lambda}
\label{5.17}
\eeq
(the sum is over $k\in\Zset_{\ge 0}$ and only a finite number of $a_k\ne0$).
For Bessel wave functions this definition does not make sense when $\lambda=0$
(since $\Psi_\beta(x,z)$ has a singularity at $z=0$ for $N>1$). In this case we
say that
$c$ is supported at $z=0$ iff it is of the form (cf. \leref{5.1} (ii, iii))
$$
\langle c,\Psi_\beta(x,z)\rangle =
\sum_{\alpha} \sum_{j} b_{\alpha j} x^{\alpha} (\ln x)^j.
$$
The sums are over $\alpha\in\bigcup_{i=1}^N \{\beta_i+N\Zset_{\ge 0}\}$ and
$0\le j \le\mult(\alpha)-1$ where $\mult(\alpha)$ is the multiplicity of
$\alpha$ in the above union (only a finite number of $b_{\alpha j} \ne0$).
The space of conditions $C$ is called {\em homogeneous\/}
iff it has a basis of homogeneous conditions $c$ (i.e.\ the support of $c$ is a
point).

It is easy to see that if $C$ is homogeneous then the spectral curve $\spec
A_W$ is {\em rational and unicursal\/} \cite{W} (i.e.\ its singularities can
be only cusps) -- the condition $c$ supported at $\lambda$ ``makes'' a cusp at
$\lambda$. For rank one algebras rationality and unicursality of $\spec A_W$
are necessary and sufficient for bispectrality \cite{W}. For rank $N>1$
another necessary condition is that $\spec A_W$ be $\Zset_N$-{\em invariant\/},
i.e.
\beq
A_W\subset \Cset[z^N].
\label{5.18}
\eeq
When $W$ is a Darboux transformation of a Bessel plane $V_\beta$, with generic
$\beta\in \Cset^N$, this condition is satisfied because of
Propositions \ref{p5.4}, \ref{p3.5}. It is natural to demand that
the space of conditions $C$ (or equivalently $\ker P$) also be
$\Zset_N$-invariant.

The $\Zset_N$-invariance of $\ker P$ simply means that
\beq
f(x)\in\ker P\;\Rightarrow\; f(\varepsilon x)\in \ker P,
\quad \varepsilon=e^{2\pi i/N}.
\label{5.19}
\eeq
It is easy to see that $C$ is {\em{homogeneous and
$\Zset_N$-invariant\/}} iff  $\ker P$
has a basis which is a union of:

(i) Several groups of elements supported at $0$ of the form:
\beq
\partial_y^l \Bigl( \sum_{k=0}^{k_0} \sum_{j=0}^{\mult(\beta_i + kN) - 1}
b_{k j} x^{\beta_i + kN} y^j \Bigr) \Big|_{y = \ln x},
\quad 0 \le l \le j_0,
\label{5.20a}
\eeq
where $j_0 = \max\{j|b_{kj}\not=0 \textrm{ for some } k\}$;

(ii) Several groups of elements supported at the points $\varepsilon^i \lambda$
($0 \le i \le N-1$, $\lambda\not=0$) of the form:
\beq
\sum_{k=0}^{k_0} a_k \varepsilon^{ki}
\partial_z^k\Psi_\beta(x,z)|_{z=\varepsilon^i\lambda},\quad 0\le i\le N-1.
\label{5.20b}
\eeq
Instead of \eqref{5.20b} we can also take
\beq
\sum_{k=0}^{k_0} a_k D^k_z\Psi_\beta(x,z)|_{z=\varepsilon^i\lambda},
 \quad 0\le i\le
N-1.
\label{5.20c}
\eeq
Denote by $n_0$ the number of conditions $c$ supported at $0$ (i.e.\ the number
of elements of the form \eqref{5.20a} in the above basis of $\ker P$). For
$1\le j\le
r$ denote by $n_j$ the number of conditions $c$ supported at each of the points
$\varepsilon^i\lambda_j$, $0\le i\le N-1$ (i.e.\ the number of groups of
elements of the form \eqref{5.20b} with $\lambda=\lambda_j$).

We have at last arrived at our fundamental definition.
\bde{5.5}
We say that the wave function $\Psi_W(x,z)$ is a
{\em polynomial Darboux transformation\/} of the Bessel wave function
$\Psi_\beta(x,z)$, $\beta\in\Cset^N$, iff \eqref{3.8} holds (for $V=V_\beta$)
with $P(x,\partial_x)$ and $g(z)$ satisfying:

(i) The corresponding space of conditions $C$ is homogeneous and
$\Zset_N$-invariant, or equivalently $\ker P$ has a basis of the form
(\ref{5.20a}, \ref{5.20b}).

(ii) The polynomial $g(z)$ is given by
\beq
g(z)=z^{n_0}\left(z^N-\lambda_1^N\right)^{n_1}\cdots
\left(z^N-\lambda_r^N\right)^{n_r}
\label{5.21}
\eeq
where $n_j$ are the numbers defined above.

We denote the set of planes $W$ satisfying (i), (ii) by $Gr_B(\beta)$ and put
$Gr^{(N)}_B=\bigcup_\beta Gr_B(\beta)$, $\beta\in\Cset^N$-generic.
\ede
We point out that the form \eqref{5.21} of $g(z)$ was introduced for $N=1$ by
Wilson \cite{W}. (Note that $g(z)=z^{n_0}\prod_{j=1}^r\prod_{i=0}^{N-1}
(z-\varepsilon^i\lambda_j)^{n_j}$.) We make this normalization in order that
$\Psi_{bW}(x,z)=\Psi_W(z,x)$ be a wave function; for the bispectral problem it
is inessential.
\bde{5.6}
We say that the polynomial Darboux transformation $\Psi_W(x,z)$ of
$\Psi_\beta(x,z)$ is {\em monomial\/} iff
$$
g(z)=z^{n_0}
$$
(i.e.\ iff all conditions $c$ are supported at 0).
Denote the set of the corresponding planes $W$ by $Gr_{MB}(\beta)$ and put
$Gr^{(N)}_{MB}=\bigcup_\beta Gr_{MB}(\beta)$, $\beta\in\Cset^N$-generic.
\ede
The next theorem provides another equivalent definition of $Gr_B(\beta)$ and
is used essentially in the proof of the bispectrality in the next section.
\bth{5.7}
 The wave function $\Psi_W(x,z)$ is a polynomial Darboux
transformation of the Bessel wave function $\Psi_\beta(x,z)$,
$\beta\in\Cset^N$, iff {\em(\ref{3.8}, \ref{3.9}, \ref{3.15},
\ref{3.16})}  hold {\em(}for $V=V_\beta${\em)} and

{\em(i)}  The operator $P$ has the form
\beq
P(x,\partial_x)=x^{-n}\sum_{k=0}^n p_k(x^N)(x\partial_x)^k,
\label{5.22}
\eeq
where $p_k$ are rational functions, $p_n\equiv 1$.

{\em(ii)}  There exists the formal limit
\beq
\lim_{x\to\infty} e^{-xz}\Psi_W(x,z)=1.
\label{5.23}
\eeq
\eth
\noindent
The {\it proof\/} will be split into three lemmas.
Before giving it we shall make a few comments.

The rationality of $P$ is always necessary for bispectrality \cite{DG, W},
\eqref{5.22} also imposes the $\Zset_N$-invariance. The condition \eqref{5.23}
is necessary in order that $\Psi_{bW}(x,z)=\Psi_W(z,x)$ be a wave function. The
limit in \eqref{5.23} is formal in the sense that it is taken in the
coefficient at any power of $z$ in the formal expansion \eqref{1.13}
separately, i.e. \beq
\lim_{x\to\infty} a_j(x)=0\quad {\rm for\ all}\ \ j\ge 1.
\label{5.24}
\eeq
Our first lemma is similar to Proposition 5.1 ((i) $\Rightarrow$
(ii)) from \cite{W}.
\ble{5.8}
 If $P$ has rational coefficients and is $\Zset_N$-invariant {\em(}see
\eqref{5.22}{\em)}
then the conditions $C$ are homogeneous and $\Zset_N$-invariant {\em(}see
{\em(\ref{5.20a}, \ref{5.20b}))}
\ele
\noindent
{\it Proof}. If $\ker P=\span \{f_0,\ldots, f_{n-1}\}$, the second coefficient
of $P$ is
$$
-\partial_x\log \Wr(f_0,\ldots , f_{n-1})
$$
and is rational. \leref{5.1} implies that $\Wr(f_0,\ldots,f_{n-1})$ is of the
form $$
x^\alpha e^{\lambda x}\times(\mbox{Laurent series in } x^{-1}).
$$
In particular each element of $\ker P$ is a sum of terms of the form
$$
 e^{\lambda x}\times(\mbox{Laurent series in } x^{-1})
\quad {\rm or} \quad
x^\alpha (\ln x)^k.
$$
We order the (finite) set of all such
$e^{\lambda x}$ and $x^\alpha(\ln x)^k$ occuring in $\ker P$. The highest term
in $\Wr(f_i)$ is just the Wronskian of the highest terms of the $f_i$. If it
vanishes
then the highest terms of the $f_i$ are linearly dependent, so by a linear
combination we can obtain a new basis with lower highest terms. So we can
suppose that the highest term of $\Wr(f_i)$ is non-zero. Repeating the same
argument with the lowest term, we shall finally obtain a basis whose elements
consist of only one term, i.e.\ are homogeneous (cf.\ \cite{W}).

Because the coefficients of $P$ are rational, \eqref{3.1} implies that it does
not matter which branch of the functions 
$x^\alpha (\ln x)^k$ in $\ker P$ we take
for $x\in\Cset$. Let
$$
\sum_{j=0}^{j_0} f_j(x) (\ln x)^j \in\ker P
$$
with
$
f_j(x) = \sum_\alpha b_{\alpha j} x^\alpha.
$
Then
$
\sum f_j(x) (\ln x + 2l\pi i)^j \in\ker P
$
for arbitrary $l\in\Zset$ and also for $\l\in\Cset$ since it is polynomial in
$l$. Taking the derivative with respect to $l$ we obtain that
$$
\sum_{j=0}^{j_0} f_j(x) j (\ln x)^{j-1}
$$
also belongs to $\ker P$.

On the other hand the $\Zset_N$-invariance of $P$ (see \eqref{5.19}) implies
$
\sum f_j(\varepsilon x) (\ln x + 2\pi i/N)^j \in\ker P
$
and
$$
\sum_{j=0}^{j_0} f_j(\varepsilon x) (\ln x)^j \in\ker P
$$
for $\varepsilon = e^{2\pi i/N}$.

Now it is obvious that $\ker P$ has a basis of the form (\ref{5.20a},
\ref{5.20b}).
\qed
\ble{5.9}
 If\/ $\ker P$ has a basis of the form {\em(\ref{5.20a}, \ref{5.20b})} then
$P$ has rational coefficients and is $\Zset_N$-invariant {\em(}see
\eqref{5.22}{\em)}.
\ele
\noindent
{\it Proof}. Consider first the case when the basis of $\ker P$ is
$$
f_i(x)=\sum_k a_k\partial_z^k\Psi_\beta (\varepsilon ^i x,z)|_{z=\lambda},
\quad 0\le i\le N-1, \; \lambda\not=0.
$$
We shall show that
$\det\left(\partial_x^{n_j} f_i(x)\right)_{0\le i, j\le N-1}
$
is a rational function of $x$ for arbitrary $n_j\in \Zset_{\ge0}$. Using
(\ref{2.17'1}, \ref{2.17'2})
we can express all derivatives
of $\Psi_\beta(x,z)$ (both with respect to $z$ and $x$) only by
$\partial_x^k\Psi_\beta(x,z)$, $0\le k\le N-1$, to obtain
\beq
\partial_x^{n_j} f_i(x) =
\sum_{k=0}^{N-1} \alpha_{kj} (x,\lambda) \partial_x^k\Psi_\beta(\varepsilon^i
x,\lambda)
\label{5.25}
\eeq
with rational coefficients $\alpha_{kj}$. Therefore
$$
\det \left(\partial_x^{n_j}f_i(x)\right)=
\det\bigl(\alpha_{kj}(x,\lambda)\bigr)
\det\bigl(\partial_x^k\Psi_\beta(\varepsilon^i x,\lambda)\bigr).
$$
But
$\det\bigl(\partial_x^k\Psi_\beta(\varepsilon^i x,\lambda)\bigr)=\const$
because the second coefficient of $L_\beta-\lambda$ is 0.
If the basis $f_0,\ldots,f_{mN-1}$ of $\ker P$ contains $m$ groups of the type
considered above (i.e.\ \eqref{5.20b}) we can represent the matrix
$$
\left(\partial_x^{n_j}f_i(x)\right)_{0\le i, j\le mN-1}, \quad n_j\in
\Zset_{\ge0},
$$
in the block-diagonal form
$$
\pmatrix{W_1&0\cr
0&W_2&\vdots\cr
&\ldots\cr
&&W_m\cr}
$$
where each block $W_s$ has the form already considered above. This can be
achieved by columns and rows operations, using the representation \eqref{5.25}.

If in addition there are some groups of elements of the form
\eqref{5.20a}, we kill the logarithms by columns operations and then cancel the
powers $x^{\beta_i}$ from the numerator and the denominator of \eqref{3.1}.
\qed
\ble{5.10}
 If $C$ is homogeneous and $h$, $g$ are as in {\em \eqref{5.1}},
{\em\eqref{5.21}}, then {\em \eqref{5.23}} is satisfied. Conversely, {\em
\eqref{5.23}} implies {\em\eqref{5.21}}.
\ele
\noindent
{\it Proof}. The second part of the lemma is an obvious consequence of the
first one.

For a basis $\{\Phi_i(x)\}_{0\le i \le dN-1}$ of $\ker h(L_\beta)$ ($d=\deg h$)
we consider the basis of $\ker P$
\beq
f_k(x)=\sum_{i=0}^{dN-1} a_{ki}\Phi_i(x),\quad 0\le k\le n-1.
\label{3.24}
\eeq
Formulae (\ref{3.8}, \ref{3.1}) imply
\beqa
\Psi_W(x,z)&=&\frac{\Wr\bigl( f_0(x),\ldots, f_{n-1}(x),\Psi_\beta(x,z)\bigr)}
{g(z)\Wr\bigl( f_0(x),\ldots, f_{n-1}(x)\bigr)}
\label{3.25} \\
&=&\frac{\sum\det A^I \Wr\bigl(\Phi_I(x)\bigr) \Psi_I(x,z)}
{\sum\det A^I \Wr\bigl(\Phi_I(x)\bigr)}.
\hfill \label{3.26}
\eeqa
The sum is  taken over all $n$-element subsets
$$
I=\{i_0<i_1<\ldots <i_{n-1}\} \subset\{0,1,\ldots, dN-1\}
$$
and here and further we use the following notation: $A$
is the matrix from \eqref{3.24} and
$
A^I=(a_{k,i_l})_{0\le k,\; l\le n-1}
$
is the corresponding minor of $A$,
$
\Phi_I(x)=\left\{\Phi_{i_0}(x),\ldots,\Phi_{i_{n-1}}(x)\right\}
$
is the corresponding subset of the basis $\{\Phi_i(x)\}$ of $\ker h(L_\beta)$
and
\beq
\Psi_I(x,z)=\frac{\Wr\bigl(\Phi_I(x),\Psi_\beta(x,z)\bigr)}{g(z)\Wr\bigl(
\Phi_I(x) \bigr)}
\label{3.27}
\eeq
is a Darboux transformation of $\Psi_\beta(x,z)$ with a basis of  $\ker P$
$f_k=\Phi_{i_k}$.

Using \eqref{3.26} it is sufficient to prove \eqref{5.23} for
$\Psi_I(x,z)$, hence we can take $\ker P$ consisting of functions
\beqa
&&f_i(x)=\partial_z^{k_i} \Psi_\beta(x,z)\big|_{z=\lambda_i},\quad 0\le i\le
p-1, \hfill\label{5.26a}\\
&&f_i(x)=x^{\alpha_i}(\ln x)^{l_i},  \quad p\le i\le n-1.
\hfill\label{5.26b}
\eeqa
We shall consider the case when $\lambda_i \neq \lambda_j$ for $i \neq j$.
The general case can be reduced to this by taking a limit.
In the formula \eqref{3.25} we expand the determinants in the last $n-p$
columns (using the Laplace rule):
\beq
\Wr(f,\Psi_\beta)= \sum \pm \det
\left(\partial_x^{j_s}f_i(x),
\partial_x^{j_s}\Psi_\beta(x,z)\right)_{{
0\le s\le p \atop 0\le i\le p-1}}
\cdot
\det\left(\partial_x^{j_s} f_i(x)\right) _{p+1\le s\le n \atop p\le i\le n-1};
\hfill\label{5.27}
\eeq
\beq
\Wr(f)=\sum\pm\det
\left(\partial_x^{j_s} f_i(x)\right)_{0\le s,i\le p-1}.
\det\left(\partial_x^{j_s}f_i(x)\right)_{p\le s, i\le n-1},
\hfill\label{5.28}
\eeq
where the sums are over the permutations $(j_0,\ldots, j_n)$ (resp.\
$(j_0,\ldots, j_{n-1})$) of $(0,\ldots,n)$ (resp.\ $(0,\ldots,n-1)$) such that
$j_0<\ldots< j_p$ and $j_{p+1}<\ldots <j_n$
(resp.\ $j_0<\ldots<j_{p-1}$ and $j_p<\ldots<j_{n-1}$). We extract the terms
with the highest power of $x$ in the numerator and in the denominator of
\eqref{3.25}. Obviously,
\beq
\det\left(\partial_x^{j_s} f_i(x)\right)_{p\le i\le n-1} =\const \cdot x^{
\sum_{i=p}^{n-1}\alpha_i -\sum_s j_s} R_J(\ln x)
\label{5.29}
\eeq
for some polynomials $R_J(\ln x)$ ($J$ is the permutation $(j_s)$).
On the other hand for $0\le i\le p-1$
\beq
\partial_x^{j_s} f_i(x)= \partial_x^{j_s} \partial_z^{k_i}
\Psi_\beta(x,z) \big|_{z=\lambda_i}
= x^{k_i} e^{\lambda_i x}\left(\lambda_i^{j_s} +O(x^{-1})\right)
\label{5.30}
\eeq
and
\beq
\partial_x^{j_s}\Psi_\beta(x,z) = e^{xz}\left(z^{j_s}+O(x^{-1})\right).
\label{5.31}
\eeq
Now it is easy to see that the  leading terms are obtained for the permutations
$$
(n-p,n-p+1,\ldots , n, 0,1,\ldots, n-p-1),
$$
respectively
$$
(n-p,n-p+1,\ldots , n-1, 0,1,\ldots, n-p-1).
$$
Substituting (\ref{5.29}, \ref{5.30}, \ref{5.31}) in (\ref{5.27}, \ref{5.28})
and canceling the determinant \eqref{5.29} for
$J = (0,1,\ldots,n-p-1)$,
we derive that
$$ \lim_{x \to\infty} e^{-xz} P(x,\partial_x) \Psi_\beta(x,z)$$
is a fraction of two van der Monde determinants and therefore is equal to
$g(z)$.
\qed
\sectionnew{Bispectrality of polynomial Darboux transformations}
In this section we prove the main result of the paper, \thref{6.3},
claiming that polynomial Darboux transformations (see \deref{5.5}),
performed on Bessel operators, produce bispectral operators.
On its hand \thref{6.3} is an almost obvious consequence of
\thref{6.2} in which we prove that the bispectral involution 
is well-defined on the submanifolds $\GrBb$ and maps them
into themselves. The importance of \thref{6.2} is not only to provide
a proof of our main result (\thref{6.3}) but also to enlighten
the bispectral involution. Its proof uses only the definition
of polynomial Darboux transformation from \thref{5.7}
(i.e.\ it does not use \deref{5.5}). On the other hand, the 
proof is completely constructive and together with \deref{5.5}
it provides an algorithmic procedure to compute bispectral
wave functions and the corresponding bispectral operators.
This procedure is described at the end of the section.
Many examples computed by making use of it are presented in Sect.~5.

Let $V_\beta$ be a Bessel plane for a generic $\beta\in\Cset^N$ (i.e.\
$V_\beta$ is not a Darboux transformation of $V_{\beta'}$ with $\beta'\in
\Cset^{N'}$, $N'<N$). In this section $W$ will be a polynomial Darboux
transformation of $V_\beta$, i.e.
$$
W\in Gr_B(\beta).
$$
We use the notation from (\ref{3.8}, \ref{3.9}) with $V=V_\beta$.
In the next proposition we show that the manifold of polynomial Darboux
transformations is preserved by the involutions $a$ and $s$ (introduced in
Subsect.~1.4).
\bpr{6.1}
If $W\in Gr_B(\beta)$, then

{\rm{(i)}} $sW\in Gr_B(\beta)$;

{\rm{(ii)}} $aW\in Gr_B(a(\beta))$, where $a(\beta)=(N-1)\delta-\beta$,
$\delta=(1,1,\ldots,1)$.
\epr
\noindent
{\it Proof}. First recall that (\prref{4.4}) $sV_\beta= V_\beta$ and
$aV_\beta=V_{a(\beta)}$. We shall study the action of the involutions on
$\Psi_W(x,z)$ and check that the conditions of \thref{5.7} are satisfied.

{\rm{(i)}} is trivial because
$$
\Psi_{sW}(x,z)=\Psi_W(-x,-z)=\frac1{g(-z)}P(-x,-\partial_x)\Psi_\beta(x,z).
$$
To prove {\rm{(ii)}} we note that the $\Zset_N$-homogeneity of $P$ (see
\eqref{5.22}) is equivalent to
\beq
P(\varepsilon x,\varepsilon^{-1}\partial_x) =\varepsilon^{-n} P(x,\partial_x),
\label{6.3}
\eeq
for $n=\ord P$, $\varepsilon=e^{2\pi i/N}$.
It follows from \eqref{3.16} that the operator $Q$ (from \eqref{3.9}) has the
same property
and also that $Q=h(L_\beta)P^{-1}$ has rational coefficients. \prref{4.3}
implies that $\Psi_{aW}$ is a Darboux transformation of $\Psi_{a(\beta)}$ with
$$
\Psi_{aW}(x,z)=\frac1{\check g(z)} Q^* (x,\partial_x) \Psi_{a(\beta)}(x,z).
$$
Obviously, $Q^*$ also satisfies \eqref{6.3}. To check \eqref{5.23}, we set
$$
K_W=1+\sum_{j=1}^\infty a_j(x)\partial^{-j}_x
$$
(see (\ref{1.13}, \ref{1.14})).
Recalling that
$$
K^*_W=1+\sum_{j=1}^\infty (-\partial_x)^{-j} a_j(x)
$$
and
$$
K_{aW}=1+\sum_{j=1}^\infty b_j(x) \partial_x^{-j}=(K^*_W)^{-1},
$$
we compute the coefficients $b_j(x)$ inductively and find that all of them are
polynomials in $a_j(x)$ and their derivatives. But by \thref{5.7} all $a_j(x)$
are {\em{rational\/}} functions of $x$ and $\lim_{x\to\infty} a_j(x)=0$, which
leads
to $\lim_{x\to\infty} b_j(x)=0$ for all $j\ge1$. This proves  \eqref{5.23} for
$aW$ (cf.\ \eqref{5.24}).
\qed

\smallskip

\prref{6.1} shows that the involutions $a$ and $s$ preserve $Gr_B^{(N)}$.
The central result of the present paper is that the bispectral involution $b$
has the same property. It immediately implies that wave functions $\Psi_W$ with
$W\in Gr_B^{(N)}$ give solutions to the bispectral problem. Our next theorem
addresses this issue.
\bth{6.2}
 If $W\in Gr_B(\beta)$ then $bW$ exists and $bW\in
Gr_B(\beta)$.
\eth
\noindent
{\it Proof.} Before proving the existence of $bW$, we shall find an analog of
\eqref{3.8} for $\Psi_{bW}(x,z) = \Psi_W(z,x)$, i.e.\ we shall show the
existence
of an operator $P_{\rm b}(x,\partial_x)$ and a polynomial $g_{\rm b}(z)$ such
that
\beq
\Psi_{bW}(x,z) =\frac{1}{g_{\rm b}(z)} P_{\rm b}(x,\partial_x)\Psi_\beta(x,z).
\label{6.4}
\eeq
{}From \eqref{5.22} it follows that the operator $P$ can be written as
\beq
P(x,\partial_x)=\frac{1}{x^n\overline p_n(x^N)}
\sum_{k=0}^n \overline p_k(x^N)(x\partial_x)^k,
\label{6.5}
\eeq
where now $\overline p_k(x^N)$ are polynomials. Use
(\ref{2.17'1}--\ref{2.17'3}) to obtain
\beqa
\Psi_W(x,z) &=& \frac{1}{x^n \overline p_n(x^N) g(z)}
\sum \overline p_k(x^N)(x\partial_x)^k \Psi_\beta(x,z) \nn\\
&=& \frac{1}{x^n \overline p_n(x^N) g(z)}
\sum(z\partial_z)^k \overline p_k\bigl(L_\beta(z,\partial_z)\bigr)
\Psi_\beta(x,z). \nn
\eeqa
This implies \eqref{6.4} with
\beqa
&&P_{\rm b}(x,\partial_x) =\frac1{g(x)}
\sum_{k=0}^n (x\partial_x)^k \overline p_k\bigl(L_\beta(x,\partial_x)\bigr),
\label{6.6}
\hfill\\
&&g_{\rm b}(z)=z^n\overline p_n(z^N).
\label{6.7} \hfill
\eeqa

Now we can prove the existence of $bW$, i.e.\ that $\Psi_{bW}(x,z)$ is a wave
function (see \eqref{1.13}). Indeed, using \eqref{6.4} we can differentiate the
formal expansion \eqref{2.14'5} of $\Psi_\beta(x,z)=\Psi_\beta(xz)$; expanding
$g_{\rm b}^{-1}(z)$ at $z=\infty$ we obtain
$$
\Psi_{bW}(x,z)=e^{xz}\sum_{k\ge k_0} b_k(x)z^{-k}
$$
for some finite $k_0$. Note that the coefficients $b_k(x)$ are rational. On the
other hand
$$
\Psi_{bW}(x,z)=\Psi_W(z,x)=e^{xz} \sum_{j\ge 0}a_j(z) x^{-j}
$$
with rational $a_j(z)$ such that (see \eqref{5.24})
\beq
\lim_{z\to\infty} a_j(z)=0,
\ j\ge1;\quad a_0(z)\equiv1.
\label{6.8}
\eeq
These two (formal) expansions  of $\Psi_{bW}(x,z)$ are connected by
$$
a_j(z)=\sum_{k\ge k_0} b_{kj} z^{-k},
$$
where
$$
b_k(x)=\sum_j b_{kj} x^{-j}, \quad b_{kj}=0 \ \ {\rm for}\ \ j<0.
$$
Now \eqref{6.8} implies
$$
b_{kj}=0\ \ {\rm for}\ \  k<0\; j\ge 1.
$$
This shows that
$$
\Psi_{bW}(x,z)=e^{xz} \Bigl( 1+\sum_{k\ge1} b_k(x)z^{-k}\Bigr)
$$
is a wave function. It is clear that it satisfies \eqref{5.23} as well.

To show an analog of \eqref{3.9}, i.e.\ that
\beq
\Psi_\beta(x,z)=\frac{1}{f_{\rm b}(z)} Q_{\rm b}(x,\partial_x)\Psi_{bW}(x,z)
\label{6.9}
\eeq
with an operator $Q_{\rm b}$ and a polynomial $f_{\rm b}$, we shall use the
above proven identity \eqref{6.4} with $asW$ instead of $W$. It follows from
\prref{4.3} that
\beq
\Psi_{asW}(x,z)= \frac{1}{f(z)} Q^*(-x,-\partial_x)\Psi_{a(\beta)}(x,z).
\label{6.10}
\eeq
\prref{6.1} and \thref{5.7} (i) allow us to present $Q^*(-x,-\partial_x)$
in the form
\beq
Q^*(-x,-\partial_x)=\frac{1}{x^m\overline q_m(x^N)}
\sum_{s=0}^m \overline q_s(x^N)(x\partial_x)^s
\label{6.11}
\eeq
with polynomials ${\overline q}_s(x^N)$. Then
\beq
\Psi_{basW}(x,z) =\frac{1}{f(x) z^m\overline q_m(z^N)}
\sum_{s=0}^m (x\partial_x)^s \overline q_s
\bigl(L_{a(\beta)}(x,\partial_x)\bigr) \Psi_{a(\beta)}(x,z).
\label{6.12}
\eeq
The identity $ab=bas$ \cite{W} and \prref{4.3} now lead to \eqref{6.9}
with
\begin{eqnarray}
Q_{\rm b}(x,\partial_x) &=& \left( \frac{1}{f(x)}
\sum_{s=0}^m (x\partial_x)^s \overline q_s \bigl(
L_{a(\beta)} (x,\partial_x) \bigr)\right)^*
\nn \\
&=& \sum_{s=0}^m \overline q_s \left(
(-1)^N L_\beta(x,\partial_x) \right) (-x\partial_x -1)^s \frac{1}{f(x)}
\label{6.13}
\end{eqnarray}
and
\beq
f_{\rm b}(z)=(-z)^m\overline q_m\left((-z)^N\right).
\label{6.14}
\eeq
{}From \eqref{5.21} and \eqref{6.6} it is obvious that $P_{\rm b}$ is
$\Zset_N$-homogeneous. This completes the proof of \thref{6.2}.
\qed

\smallskip

An immediate corollary is the following result, which we state as a theorem
because of its fundamental character.
\bth{6.3}
If  $W\in Gr^{(N)}_B$ then the wave function
$\Psi_W(x,z)$ solves the bispectral problem, i.e.\ there exist operators
$L(x,\partial_x)$ and $\Lambda(z,\partial_z)$ such that
\beqa
&&L(x,\partial_x)\Psi_W(x,z) = h(z^N)\Psi_W(x,z),
\label{6.15}\\
&&\Lambda(z,\partial_z)\Psi_W(x,z) = \Theta(x^N)\Psi_W(x,z),
\label{6.16}
\eeqa
Moreover,
\beq
\rank \A_W=\rank \A_{b W}=N.
\label{6.17}
\eeq
\eth
\noindent
{\it Proof.} (\ref{6.15}, \ref{6.16}) follow from
(\ref{3.8}, \ref{3.9}, \ref{6.4}, \ref{6.9}) if we set
\beqa
&&L(x,\partial_x)=P(x,\partial_x)Q(x,\partial_x),
\quad h(z^N)=f(z)g(z);
\label{6.18} \\
&&\Lambda(z,\partial_z)=P_{\rm b}(z,\partial_z)Q_{\rm b}(z,\partial_z),
\quad \Theta(x^N)=f_{\rm b}(x)g_{\rm b}(x).
\label{6.19}
\eeqa
The eq. \eqref{6.17} follows from Propositions \ref{p3.5} (i) and \ref{p5.4}.
\qed
\bex{9.7}
All {\it bispectral algebras of rank 1\/} are polynomial
Darboux transformations of the plane $H_+=\{z^k\}_{k \geq 0}$ (see \cite{W}).
This corresponds to the $N=1$ Bessel with
$$
\beta=(0),\quad L_{(0)}=\partial_x, \quad V_{(0)}=H_+=\{z^k\}_{k \geq 0}, \quad
\psi_{(0)}(x,z)=e^{xz}.
$$
Every linear functional on $H_+$ is a linear combination of
$$
e(k,\lambda)=\partial_z^k|_{z=\lambda}
$$
and $h\left(L_{(0)}\right) =h(\partial_x)$ is an operator with constant
coefficients. The ``adelic Grassmannian'' $Gr^{ad}$, introduced by Wilson
\cite{W}, coincides with $Gr_B((0))$ $(=Gr_B^{(1)})$.
In our terminology the result of \cite{W} can be reformulated as follows.

{\em{All bispectral operators belonging to rank one bispectral algebras are
polynomial Darboux transformations of operators with constant coefficients}}.
\qed
\eex
\bre{conv}
The eigenfunction $\Psi_W(x,z)$ from eq. \eqref{3.8} is a formal series. Let
$\Phi_\beta(x,z) = \Phi_\beta(xz)$, where $\Phi_\beta(z)$ is the Meijer's
$G$-function \eqref{1.33'2} (or any convergent solution of \eqref{2.14} in
arbitrary domain) and set
\beq
\Phi_W(x,z)=\frac{1}{g(z)} P(x,\partial_x) \Phi_\beta(x,z).
\label{3.18a}
\eeq
Then
\beq
\Phi_W(x,z)=\frac{1}{g_{\rm{b}}(x)} P_{\rm{b}}(z,\partial_z) \Phi_\beta(x,z)
\label{3.18b}
\eeq
because of \eqref{2.14} and $x \partial_x \Phi_\beta(x,z)=
z \partial_z \Phi_\beta(x,z)$. The equations $Q P=h(L_\beta)$ and
$Q_{\rm{b}} P_{\rm{b}}=\Theta(L_\beta)$ imply
\beqa
&&\Phi_\beta(x,z)=\frac{1}{f(z)} Q(x,\partial_x) \Phi_\beta(x,z),
\label{3.18c}\\
&&\Phi_\beta(x,z)=\frac{1}{f_{\rm{b}}(x)} Q_{\rm{b}}(z,\partial_z)
\Phi_W(x,z).
\label{3.18d}
\eeqa
So, we proved  that $\Phi_W(x,z)$ is a convergent bispectral eigenfunction of
the same operators $L(x, \partial_x)$ and $\Lambda(z, \partial_z)$ as
$\Psi_W(x,z).$ The involutions  $a$, $s$ and $b$ can be defined on the
manifold of ``convergent'' polynomial Darboux transformations \eqref{3.18a}
by the equations (\ref{4.7}, \ref{4.8}, \ref{1.42'5}) in which $\Psi$ is
replaced by $\Phi$ and they preserve it
(\prref{4.3} (i) now becomes a definition).
The validity of the equation $ab=bas$ in the ``convergent'' case is a
consequence of that in the ``formal'' one (see the proof of \thref{6.2}).
The rationality of the coefficients of the operator $P(x, \partial_x)$ implies
that its kernel has one and the same form
(see eqs.\ (\ref{5.20a}, \ref{5.20b})) in $\Psi$- and in $\Phi$-bases.
\qed
\ere
It is not difficult to provide an explicit algorithm for producing bispectral
pairs $L(x,\partial_x)$, $\Lambda(z,\partial_z)$. Although obvious we have
collected the steps of this algorithm as they are scattered in the present and
the previous sections.

\medskip
\noindent
{\it Step 1}. Choose an arbitrary set of conditions based in some points
$\lambda_0=0$, $\lambda_1,\ldots,\lambda_r$ of the form
(\ref{5.20a}, \ref{5.20b}), i.e.
a basis of $\ker P$. The proof of \leref{5.9} provides an explicit computation
of
the coefficients of $P$ in terms of the coefficients $a_k$, $b_{kj}$ in $\ker
P$. The polynomial $g(z)$ is given by \deref{5.5} (ii).

\medskip
\noindent
{\it Step 2}. Take $h(z^N)= z^{d_0N} \prod_{j=1}^r\left(z^N-\lambda_j^N\right)
^{d_j}$ with high enough powers $d_0,\ldots,d_r$ such that $\ker P \subset\ker
h(L_\beta)$ (cf.\ \leref{5.1}).
The minimal such $d_j$'s can be computed as follows.

(i) For a condition, supported at $0$, of the form \eqref{5.20a} set $j(k) =
\max\{j | b_{kj}\not=0 \}$, $0\le k\le k_0$. Let $\beta_i + kN = \beta_{i_s} +
p_s N$ for $0\le s\le \mult(\beta_i + kN)-1$ with $0\le p_0\le\ldots\le
p_{\mult(\beta_i+kN)-1}$ and $i_s \not= i_t$ for $s\not= t$. Then set
$$
d_0 = 1 + \max p_{j(k)},
$$
the maximum is over all $k$ and all conditions of the form \eqref{5.20a}.

(ii) For a condition, supported at $\lambda_j\not=0$, of the form \eqref{5.20b}
let $k_0 = \max \{k | a_k\not=0 \}$. Then set
$$
d_j = 1 + \max k_0,
$$
the maximum is over all conditions of the form \eqref{5.20b} supported at
$\lambda_j$.

Then put $f(z) = h(z^N) / g(z)$.

\medskip
\noindent
{\it Step 3}. Find the coefficients of the operator $Q(x,\partial_x)$
recursively out of the equation $Q(x, \partial_x)P(x, \partial_x)=h(L_\beta(x,
\partial_x))$. Then $L(x,\partial_x)=P(x,\partial_x) Q(x,\partial_x)$.
A lower order operator $L$ can be constructed using \prref{3.5},
i.e.\ find $u(L_\beta)$ such that $\ker P$ is invariant under $u(L_\beta)$
and then $L$ out of the equation $L P = P u(L_\beta)$.

\medskip
\noindent
{\it Step 4}. Compute by \eqref{6.6} $P_{\rm b}(x,\partial_x)$ and by
\eqref{6.7}
$g_{\rm
b}(z)$. Also \eqref{6.13} and \eqref{6.14} give $Q_{\rm b}(x,\partial_x)$ and
$f_{\rm
b}(z)$. All expressions are explicit in terms of the coefficients of the
operators $P$ and $Q$. Then $\Lambda(z,\partial_z)=P_{\rm b}(z,\partial_z)
Q_{\rm b}(z,\partial_z)$ and $\Theta(x)=f_{\rm b}(x)g_{\rm b}(x)$.
\sectionnew{ Polynomial Darboux transformations of Airy planes }
This section contains analogs of the results from Sections
2 and 3 but here the building blocks are (generalized) Airy
operators (see \cite{KS, Dij}) instead of Bessel ones.
There is a minor difference in the organization of the present section
compared to that of Sections 2 and 3. Here we give the definition of
polynomial Darboux transformations on Airy wave functions
(see Definitions \ref{d8.2}, \ref{d8.3}) in the spirit of the one provided
by \thref{5.7}. Then we prove our main result \thref{8.5} (which is an 
analog of \thref{6.2}). As in Sect.~2, it automatically implies
bispectrality of the polynomial Darboux transformations.
At the end, in \prref{8.9} we show that \deref{8.3} is equivalent to
a second one (analog of \deref{5.5}) in terms of conditions on 
Airy planes. This is again important for algorithmic computations,
some of which are presented in the next section.

First we recall the definition of (generalized higher) Airy
functions. \\
For
$
\alpha = (\alpha_0,\alpha_2,\alpha_3,\ldots,\alpha_{N-1}) \in\Cset^{N-1}
$,
$\a_0 \not=0$,
consider the {\em{Airy operator}}
    \beq
L_{\alpha}(x,\partial_x) = \partial_x^N - \alpha_0 x + \sum_{i=2}^{N-1}
\alpha_i \partial_x^{N-i} \equiv P_{\alpha'}(\partial_x) - \alpha_0 x
\label{8.1}
    \eeq
where
$
\alpha' = (\alpha_2,\alpha_3,\ldots,\alpha_{N-1}).
$
 The {\em{Airy
equation\/}} is
    \beq
L_{\alpha}(x,\partial_x) \Phi(x) = 0, \text{ i.e.} \quad
P_{\alpha'}(\partial_x) \Phi(x) = \alpha_0 x \Phi(x).
\label{8.2}
    \eeq
\bex{8.1}
When $\alpha_0 = 1, \; \alpha' = 0$ eq.\ \eqref{8.2} becomes the classical
higher Airy equation (cf.\ \cite{KS})
    \beq
\partial_x^N \Phi(x) = x \Phi(x).
\label{8.3}
    \eeq
In every sector $S$ with a center at $x=\infty$  and an angle less than
$N\pi / (N+1)$, it has a solution with an asymptotics of the form
(see e.g. \cite{Wa})
     \beq
\Phi(x) \sim x^{ -\frac{N-1}{2N} } e^{ \frac{N}{N+1} x^{\frac{N+1}{N}} }
\Bigl( 1 + \sum_{i=1}^{\infty} a_i x^{-i/N} \Bigr),
\quad |x| \to \infty, \; x \in S.
\label{8.4}
\eeq
\qed
\eex
Similarly, in each sector $S$ as in \exref{8.1} eq.\ \eqref{8.2} has a solution
with an asymptotics of the form
   \beq
\Phi(x) \sim \Psi_{\alpha}(x) := x^{d/N} e^{Q(x^{1/N})}
\Bigl( 1 + \sum_{i=1}^{\infty} a_i x^{-i/N} \Bigr),
\quad |x| \to \infty, \; x \in S
\label{8.5}
   \eeq
for some $d \in\Cset$ and a polynomial $Q(x)$ of degree $N+1$ with leading
coefficient $\mu_0 \frac{N}{N+1} x^{N+1}$, where $\alpha_0 = \mu_0^N.$
The solution $\Phi$ is by no means unique, but $d$, $Q$ and all $a_i$ are
uniquely
determined and do not depend on $S$. In the sequel we shall deal only with
$\Psi_\alpha$, which is a formal solution of eq. \eqref{8.2}.
\bde{8.2}
For each $\alpha\in\Cset^{N-1}$ we call an {\em{Airy wave function\/}} the
following function
    \beq
\psi_\alpha(x,z) := \mu_0^d z^{-d} e^{-Q(\mu_0^{-1} z)} \Psi_\alpha(x,z),
\label{8.6}
    \eeq
where
$$\Psi_\alpha(x,z) := \Psi_\alpha( \alpha_0^{-1} z^N + x ).$$
\hfill \\
It is easy to see that $\psi_\alpha$ is indeed a wave function if we expand
$\Psi_\alpha( \alpha_0^{-1} z^N + x )$ at $x=0$:
    \beq
( \alpha_0^{-1} z^N + x )^{ -i/N } = \sum_{k\geq 0} {-i/N \choose k}
(\mu_0^{-1} z)^{-i-kN} x^k
\label{8.7}
    \eeq
(we shall always use $\mu_0$ as an $N$-th root of $\a_0$).

The plane in Sato's Grassmannian corresponding to $\psi_\a(x,z)$  will be
called an {\em{Airy plane\/}} and will be denoted by $V_\a$.
\ede
Obviously, $\Psi_\a(x,z)$ solves the bispectral problem
     \beqa
&&L_\a(x,\partial_x) \Psi_\a(x,z) = z^N \Psi_\a(x,z)
\label{8.8}
\hfill\\
&&L_\a( \a_0^{-1}z^N, \partial_{\a_0^{-1} z^N} )
\Psi_\a(x,z) = \a_0x \Psi_\a(x,z)
\label{8.9}
     \eeqa
because
     \beq
\partial_x \Psi_\a(x,z) = \partial_{\a_0^{-1} z^N} \Psi_\a(x,z).
\label{8.10}
     \eeq

It is clear that $\psi_\a$ satisfies \eqref{8.8} and analogs of (\ref{8.9},
\ref{8.10}) obtained by conjugating by $z^{-d} e^{-Q(\mu_0^{-1} z)}$. (Up to
this conjugation
$\Psi_\a$ and $\psi_\a$ give one and the same solution to the bispectral
problem.)

We shall define
{\em{polynomial Darboux transformations\/}} of Airy planes as in the Bessel
case (see \deref{5.5} and \thref{5.7}). Before that we shall define a
{\em{bispectral involution\/}} $b_1$ on them. Note that the involution $b$ from
\cite{W} (see Subsect.~1.4) is not well defined on $V_\a$ (i.e. $\psi_\a(z,x)$
is not a wave function). The properties of $b$ we would like $b_1$ to have,
are:

1) it has to interchange the roles of $x$ and $z$;

2) it has to preserve Airy planes. \\
Therefore we define
   \beq
b_1 \Psi_\a(x,z) := \Psi_\a(x,z) = \Psi_\a( \a_0^{-1} z^N, \mu_0 x^{1/N} ),
\label{8.11}
   \eeq
or equivalently,
   \beq
b_1 \psi_\a(x,z) := \psi_\a(x,z) = \mu_0^d x^{d/N} z^{-d} e^{ Q(\mu_0 x^{1/N})
- Q(\mu_0^{-1} z) } \psi_\a( \a_0^{-1} z^N, \mu_0 x^{1/N} ).
\label{8.12}
   \eeq
For a Darboux transformation $W$ of $V_\a$ we define $\psi_{b_1 W}$ and
$\Psi_{b_1 W}$ in a
similar way. (We still do not know whether $b_1 W \in\Gr$, the notation
$\psi_{b_1 W}$ is still  formal.)
\bde{8.3}
A Darboux transformation $W$ of 
an Airy plane $V_\a$ is called {\em{polynomial\/}}
iff (in the notation of \deref{3.2})

(i) the operator $P$ has rational coefficients;

(ii) $g(z) = g_1(z^N),$ $g_1 \in\Cset[z]$;

(iii) $\lim\limits_{z \to\infty} e^{-xz} \psi_{b_1 W}(x,z) = 1$. \\
(The limit is formal and has the same meaning as in \eqref{5.23}.)

Denote the set of all such $W\in\Gr$ by $\GrAa$ and put $\GRAN =
\bigcup_{\a\in\Cset^{N-1}} \GrAa.$
\ede
\bre{8.4}
The parts (i) and (ii) of the above definition remain the same if we
substitute $\psi_\a$ and $\psi_W$ by $\Psi_\a$ and $\Psi_W$, where
    \beq
\psi_W(x,z) := \mu_0^d z^{-d} e^{-Q(\mu_0^{-1} z)} \Psi_W(x,z).
\label{8.13}
    \eeq
\ere
The main result of this section is that $\GrAa$ is preserved by the involution
$b_1$.
\bth{8.5}
{\rm{(i)}} If $W\in\GrAa$, then $\psi_{b_1W}(x,z)$ is a wave function
corresponding to a plane $b_1 W \in\GrAa.$

{\rm{(ii)}} For $\a\in\Cset^{N-1}$ the spectral algebra $\A_{V_\a}$ is
$\Cset[L_\a].$
\eth
An immediate {\bf{corollary}} is that  {\em{the planes $W\in\GrAa$ give
solutions to the bispectral problem of\/}} $\rank N$:
$$ \rank \A_W = \rank \A_{b_1 W} = N.$$
The {\em{proof of \thref{8.5}\/}} is completely parallel to that of
\thref{6.2}. We shall be very brief, indicating only the major differences and
the most important steps. We start with a lemma illuminating the purpose of the
constraints (i) and (ii) in \deref{8.3} (cf. \eqref{6.4}).
\ble{8.6}
If $W\in\GrAa$, then
     \beq
\Psi_{b_1 W}(x,z) = \frac{1}{ g_{\mathrm b}(z) }
P_{\mathrm b}(x,\partial_x) \Psi_\a(x,z),
\label{8.14}
     \eeq
$P_{\mathrm b}$ is with rational coefficients and $g_{\mathrm b}$ is polynomial
in $z^N$.
\ele
\proof
We compute
     \beqa
&&\Psi_{b_1 W}(x,z) = \Psi_W( \a_0^{-1} z^N, \mu_0 x^{1/N} ) \nn \\
&&= \frac{ P( \a_0^{-1} z^N, \partial_{\a_0^{-1} z^N})
\Psi_\a( \a_0^{-1} z^N, \mu_0 x^{1/N} ) }{ g(\mu_0 x^{1/N}) } =
\frac{1}{g_{\mathrm b}(z)} P_{\mathrm b}(x,\partial_x)
\Psi_\a(x,z),   \nn
    \eeqa
where if
    \beq
P(x,\partial_x) = \frac{1}{ p_n(x)} \sum_{k=0}^n p_k(x) \partial_x^k,
\qquad g(z) = g_1(z^N)
\label{8.15}
    \eeq
with polynomials $p_k$ and $g_1$, then (using (\ref{8.8}, \ref{8.10}))
    \beqa
&&P_{\mathrm b}(x,\partial_x) = \frac{1}{g_1(\a_0 x)} \sum_{k=0}^n
\partial_x^k p_k(\a_0^{-1} L_\a(x,\partial_x)),
\label{8.16} \\
&&g_{\mathrm b}(z) = p_n (\a_0^{-1} z^N).
\label{8.16'5}
\hfill
\eeqa
\qed

\smallskip

The proof that $\psi_{b_1 W}(x,z)$ is a wave function is the same as in the
Bessel
case, using the above lemma and the condition (iii) of \deref{8.3}. Now the
identity  $ab=bas$ \cite{W} is modified in the following way.

Introduce the maps $p$ and $p^{-1}$ as follows
    \beqa
&&\Psi_{p W}(x,z) := \Psi_W( \a_0^{-1} x^N, \mu_0 z^{1/N} ), \nn \\
&&\Psi_{p^{-1} W}(x,z) := \Psi_W( \mu_0 x^{1/N}, \a_0^{-1} z^N ). \nn
    \eeqa
The notation $p W$, $p^{-1} W$ is formal -- these are not planes in $\Gr$. But
$$\Psi_{b_1 W}(x,z) = \Psi_{bp W}(x,z) = \Psi_{p^{-1}b W}(x,z)$$
corresponds to the wave function $\psi_{b_1 W}(x,z)$ and to $b_1 W \in\Gr$.

Multiplying the identity $ab=bas$ on the right by $p$, we obtain
    \beq
ab_1 = b_1 a_1, \quad\text{where} \;\; a_1 = p^{-1}asp.
\label{8.17}
    \eeq
Note that for $W\in\GrAa \quad aW,b_1 W$ and hence $a_1 W$ are planes in $\Gr$.

The next lemma gives the action  of the involutions on the Airy planes (the
proof is the same as that of \prref{4.4}).
\ble{8.7}
{\rm{(i)}} $sV_\a = V_{s(\a)}$, where
$s(\a) =  ( (-1)^{N+1}\a_0,\a_2,$$-\a_3,\ldots,(-1)^{N-1}\a_{N-1} );$
{\rm{(ii)}} $aV_\a = a_1 V_\a = V_{a(\a)}$, where
$a(\a) = ( (-1)^N\a_0,\a_2,-\a_3,\ldots,(-1)^{N-1}\a_{N-1} ).$
\ele
We also need an analog of \prref{6.1}.
\ble{8.8}
If $W\in\GrAa$, then $aW$ and $a_1 W$ belong to $\GrA{a(\a)}$.
\ele
For the {\em{proof\/}} we need an analog of \prref{4.3} for $a_1$. A simple
computation shows that if
  $$\Psi_V(x,z) = \frac{1}{f(z)} Q(x,\partial_x) \Psi_W(x,z)$$
for $V,W \in\GrAa$ and
     $$Q(x,\partial_x) = \sum q_k(x) \partial_x^k$$
then
  $$\Psi_{a_1 W}(x,z) = \frac{1}{f(z)}
         Q^{*_1}(x,\partial_x) \Psi_{a_1 V}(x,z)$$
with
$$Q^{*_1} = \sum \Bigl( \partial_x + \frac{1-N}{N}x^{-1} \Bigr)^k (-1)^{(N-1)k}
q_k( (-1)^N x ).$$
The rest of the proof is left to the reader.
\qed

\smallskip

The proof of part (i) of \thref{8.5} is completed exactly
as in the Bessel case.
For the part (ii), we note that while the Bessel wave functions are ``{\it
{multiplication invariant\/}}'', the Airy ones are ``{\it {translation
invariant\/}}''. More precisely, for arbitrary $ c \in \Cset $
$$\Psi_{\alpha}(x+c,(z^N-\alpha_0c)^{1/N})=
  \Psi_{\alpha}(\alpha_0^{-1}(z^N-\alpha_0c)+x+c)=
  \Psi_{\alpha}(\alpha_0^{-1}z^N+x)=
  \Psi_{\alpha}(x,z) $$
(expand $(z^N-\alpha_0c)^{1/N} = \sum_{k \geq 0}{1/N \choose k }z^{1-kN}
(-\alpha_0c)^k)$.
Let $u(z) \in A_\alpha,$ $L(x,\partial_x)\in \A_\alpha$ and
$$
L(x,\partial_x)\Psi_{\alpha}(x,z) = u(z)\Psi_{\alpha}(x,z)
$$
(this is equivalent to
$ L\psi_{\alpha}(x,z) = u \psi_{\alpha}(x,z)$). Then
$$ L(x+c,\partial_x) \Psi_{\alpha}(x,z) = u( (z^N - \alpha_0
c)^{1/N} ) \Psi_{\alpha}(x,z)$$
and
$L(x+c,\partial_x) \in \A_{\alpha}, \quad u((z^N-\alpha_0c)^{1/N})
\in A_{\alpha}.$
But $A_{\alpha} \subset \Cset[z]$, therefore $u((z^N-\alpha_0c)^{1/N})
\in \Cset[z]$ for all $c$ and $u(z) \in \Cset[z^N]$.

This completes the proof of \thref{8.5}.
\hfill \\

At the end of this section we note that an equivalent
 definition of $\GrAa$ can be given in terms of
{\it {conditions}} $C$ (cf. Sect.~2).

Using the translation invariance of $\Psi_{\a}$ we can suppose that none of
the conditions $C$ is supported at $0$. Then we have an analog of
\thref{5.7}.
\bpr{8.9}
The Darboux transformation $W$ of $V_{\a}$ is polynomial iff

{\rm {(i)}}
The space of conditions $C$ is homogeneous and
$\Zset_N$-invariant.
Equivalently, $\ker P$ has a basis of the form
\beqa
f_{ij}(x)&=& \sum_k a_{ki} \partial_{z}^k
\psi_{\a}(x,\varepsilon^j z)|_{z=\lambda_i}
\nn\\
&=&\sum_k a_{ki} \varepsilon^{-jk} \partial_z^k
\psi_{\a}(x,z) \Big| _{z=\varepsilon^j \lambda_i},
\label{8.18}
\eeqa
$ 0 \leq j \leq N-1,$ $ 1 \leq i \leq r $ {\rm (}for some $r${\rm )},
$\lambda _i \not= 0 $.

{\rm {(ii)}}
The polynomial $g(z)$ has the form \eqref{5.21},
i.e.
$$g(z) = (z^N-\lambda_1^N)^{n_1} \cdots
(z^N-\lambda_r^N)^{n_r}$$
where $n_i$ is the number of conditions $C$ supported at
each of the points $\varepsilon ^j \lambda_i, \;\; 0 \leq
j \leq N-1.$
\epr
\noindent
The {\it {proof}} of  the ``if'' part is the same as in the
Bessel case and will be omitted. (In fact, most of the proofs in
Sect.~2 are valid in a more general situation.)

The ``only if'' part is also similar to the corresponding
result in the Bessel case but some more explanation is
needed. For fixed $\lambda \not= 0$ we shall use
representations of the kernel of the operator $L_\alpha
-\lambda ^N$ in three different linear spaces of formal power
series. First set
\beq
\varphi_{\a} (x,\lambda)= \mu_0^d \lambda^{-d} e^{-Q(\lambda)}
\Psi_{\a}(y), \qquad y=\a _0^{-1} \lambda ^N+x,
\label{8.19}
\eeq
 considered as a formal power series in $y^{-1/N}$ (where
$\Psi _{\a}$ is from \eqref{8.5}). The Airy wave function
$\psi _{\a} (x,\lambda)$ (see \eqref{8.6}) is given by the same
formula after expanding $y^{-1/N}$ at $x=0$ as in
\eqref{8.7}.
The other possibility is to expand $y^{-1/N}$ at $x = \infty$:
\beq
 (x+ \a _0 ^{-1} \lambda ^N) ^{-i/N} = \sum _{k \geq 0} {-i/N
\choose k } x ^{-i/N-k} (\a_0^{-1} \lambda ^N)^k.
\label{8.20}
\eeq
Inserting \eqref{8.20} in \eqref{8.19}, we obtain another formal
series $\chi _{\a}(x,\lambda)$. Denote by $\varphi_{\a}^{(j)}$,
 $\psi_{\a}^{(j)}$, $\chi_{\a}^{(j)}$ the images of
$\varphi _{\a},$ $\psi_{\a},$ $\chi_{\a}$ under the
transformations
$$ y^{1/N} \mapsto \varepsilon ^j y^{1/N}, \quad
   \lambda \mapsto \varepsilon ^j \lambda, \quad
   x^{1/N} \mapsto \varepsilon ^j x^{1/N},$$
respectively $(\varepsilon =e^{2 \pi i/N}).$ Then $\psi_{\a}^{(j)}$
and $\chi_{\a}^{(j)}$ are obtained by expanding $\varphi_{\a}^{(j)}$ and
in the corresponding spaces of formal series
$ \ker (L_{\a}- \lambda ^N)^d$ has bases
$$ \{ \partial _{\lambda}^k \psi _{\a}^{(j)} \}, \;\;
   \{ \partial _{\lambda}^k \varphi _{\a}^{(j)} \}, \;\;
   \{ \partial _{\lambda}^k \chi _{\a}^{(j)} \}, \qquad
   0 \leq k \leq d-1, \; 0 \leq j \leq N-1. $$
Our observation is that if $ \ker P $ has a basis
$$
f_i(x) = \sum _{k,j} a_{kj}^i \partial_{z}^k
\psi_{\a}^{(j)}(x,z)|_{z=\lambda},
$$
then the same formula gives a basis of $\ker P$ when
$\psi$'s are substituted by $\varphi$'s or $\chi$'s and vice
versa. Indeed, this follows from \eqref{3.1} and the fact
that $P$ has rational coefficients.
We complete the proof of \prref{8.9} noting that
while $P$ depends rationally on $x,$ $\chi _{\a}^{(j)}$ are
formal series in $x^{1/N}$ and the same argument as in the
Bessel case gives that $\ker P$ has a $\chi$-basis of the
form \eqref{8.18}.
\qed
\sectionnew{Explicit formulae and examples}
In this section we have collected several classes of examples. We wanted at
least to include all previously known examples (unless by ignorance we miss
some of them) -- see \cite{DG, W, Z, G3, LP}. We hope that we have elucidated
and unified them.

For monomial transformations we derive formulae expressing the operators $L$
and $\Lambda$, solving the bispectral problem, only in terms of the matrix $A$
and the vector $\gamma$ (see \prref{9.1} below). This
explicit expression for $\Lambda$ (though possibly of high order) to the best
of our knowledge is new even for $N=2$ (see \cite{DG}). In other examples we
illustrate the properties of the operator of minimal order from a bispectral
algebra: when does its order coincide with the rank of the algebra and
when this operator is a Darboux transformation of a power of a Bessel operator.
We also point out that the classical Bessel potentials $u(x)=cx^{-2}$
\cite{DG} can produce new solutions of the bispectral problem for any $c$.

We describe in detail the polynomial Darboux transformations from
$(L_\a - \lambda^N)^2$ where $L_\a$ is an arbitrary Airy or Bessel operator of
order $N$. We do not want simply to show that our procedure of constructing
bispectral operators works but to point out that the
involutions $a$ and $b$ ($b_1$ in Airy case) possess some very interesting
properties which deserve further study.
\subsection{Monomial Darboux transformations of Bessel planes}
Let $\beta\in \Cset^N$ and $W\in Gr_{MB}(\beta)$. We use the notation from
(\ref{3.8}, \ref{3.9}) (with $V=V_\beta$) and from (\ref{6.4}, \ref{6.9}).
When the Darboux transformation is monomial
\beq
g(z)=z^n,\quad h(z)=z^d
\label{9.4}
\eeq
for some $n$, $d$. We shall consider only the case when there are no logarithms
in the basis \eqref{5.20a} of $\ker P$. The general case can be reduced to this
one by taking a limit in all formulae (see \exref{log} below). Now $\ker P$
has a basis of the form
\beq
f_k(x)=\sum_{i=1}^{dN} a_{ki} x^{\gamma_i},\quad 0\le k\le n-1,
\label{9.5}
\eeq
such that
\beq
\gamma_i-\gamma_j \in N\Zset\setminus 0
\qquad {\rm if}\ a_{ki}a_{kj}\not=0, \; i \not= j,
\label{7.4'5}
\eeq
where $\gamma=\beta^d$ is from \eqref{5.3}.

Let $A$ be the matrix $(a_{ki})$. We shall use multi-index notation for
subsets $I=\{i_0<\ldots< i_{n-1}\}$ of $\{1,\ldots,dN\}$ and $\delta_I$ from
\eqref{5.7}. We also put $\gamma_I = \{\gamma_i\}_{i\in I}$,
$$
A^I=(a_{k,i_l})_{0\le k, \; l\le n-1}
$$
and
$$
\Delta_I=\prod_{r<s}(\gamma_{i_r}-\gamma_{i_s}).
$$
Let $I_{\min}$ be the subset of $\{1,\ldots,dN\}$ with $n$ elements such that
$\det A^{I_{\min}}\not=0$ and $\sum_{i\in I_{\min}} \gamma_i$ be the minimum of
all such sums, and set
$$
p_I=\sum_{i\in I}\gamma_i -\sum_{i\in I_{\min}}\gamma_i.
$$
Eq. \eqref{7.4'5}  implies that these numbers are divisible by $N$. Finally,
for a subset $I$ of $\{1,\ldots, dN\}$ denote by $I^0$ its complement.

In the following proposition we express everything entering
(\ref{3.8}, \ref{3.9}, \ref{6.4}, \ref{6.9}) only in terms
of the matrix $A$ and the vector $\gamma$. Therefore for each $A$ and
$\beta\in\Cset^N$ satisfying \eqref{7.4'5} (with $\gamma=\beta^d$) we give an
explicit solution to the bispectral problem (cf. (\ref{6.18}, \ref{6.19})).
\bpr{9.1}
In the above notation the operators and the polynomials from {\rm{(\ref{6.18},
\ref{6.19})}} are given by the following formulae{\rm{:}}
\begin{eqnarray*}
&&
\hspace*{-1cm}
{\rm{(a)}} \;\;\textstyle  g(z)=z^n,
\hfill\\
&&P= \Bigl( \sum\det A^I\Delta_I x^{p_I} \Bigr)^{-1} \Bigl( \sum\det
A^I\Delta_I x^{p_I}L_{\gamma_I}\Bigr).
\hfill\\
&&
\hfill\\
&&
\hspace*{-1cm}
{\rm{(b)}} \;\;\textstyle  f(z)= z^{dN-n},
\hfill\\
&&Q= \Bigl( \sum\det A^I\Delta_I L_{\gamma_{I^0} -n\delta_{I^0}} x^{p_I} \Bigr)
\Bigl( \sum \det A^I\Delta_I x^{p_I} \Bigr)^{-1}.
\hfill\\
&&
\hfill\\
&&
\hspace*{-1cm}
{\rm{(c)}} \;\;\textstyle  g_{\rm b}(z) =z^n \sum\det A^I\Delta_I
z^{p_I},
\hfill\\
&&P_{\rm b} =\sum\det A^I\Delta_I L_{\gamma_I} (L_\beta)^{p_I/N}.
\hfill\\
&&
\hfill\\
&&
\hspace*{-1cm}
{\rm{(d)}} \;\;\textstyle f_{\rm b}(z)= z^{dN-n} \sum\det A^I\Delta_I
z^{p_I},
\hfill\\
&&Q_{\rm b}=\sum\det A^I\Delta_I (L_\beta)^{p_I/N}
L_{\gamma_{I^0}-n\delta_{I^0}}.
\end{eqnarray*}
\epr
\noindent
{\it Proof}.
Note that (c) and (d) follow from (a) and (b) (see the proof of \thref{6.2}).
To prove (a) we note that
\beqa
\Psi_W(x,z)&=&\frac{\Wr\bigl( f_0(x),\ldots, f_{n-1}(x),\Psi_\beta(x,z)\bigr)}
{z^n \Wr\bigl( f_0(x),\ldots, f_{n-1}(x)\bigr)}
\nn\\
&=&\frac{\sum\det A^I \Wr\bigl(x^{\gamma_I}\bigr) \Psi_I(x,z)}
{\sum\det A^I \Wr\bigl(x^{\gamma_I}\bigr)}.
\hfill \label{3.26'}
\eeqa
The sum is taken over all $n$-element subsets
$
I=\{i_0<i_1<\ldots <i_{n-1}\} \subset\{0,1,\ldots, dN-1\}
$,
$x^{\gamma_I} = \{x^{\gamma_i}\}_{i\in I}$ and $\Psi_I(x,z)$ are the Bessel
wave functions \eqref{5.6}. Using \eqref{5.9} and the simple fact
\beq
\Wr(x^{\gamma_I})=\Delta_I x^{\sum_{i\in I}\gamma_i -\frac{n(n-1)}2},
\label{7.16}
\eeq
we obtain (a).

To prove (b) we shall apply the involution $a$ directly on the
{\em{tau-function\/}} $\tau_W$ of the plane $W$. Recall that \cite{S,SW}
\beq
\Psi_W(t,z)= e^{\sum_{k=1}^\infty t_kz^k}
\frac{\tau\left(t-[z^{-1}]\right)}{\tau(t)},
\label{1.7}
\eeq
where $[z^{-1}]$ is the vector $\left(z^{-1}, z^{-2}/2,\ldots\right)$.
The action of $a$ is given by \cite{W}
$$
\tau_{aW}(t_1, t_2,\ldots,t_k,\ldots)=\tau_W(t_1,-t_2,
\ldots,(-1)^{k-1}t_k,\ldots).
$$
We shall need the formulae \cite{BHY2}
\beq
\tau_W(t)=\frac{\sum\det A^I\Delta_I\tau_I(t)}{\sum\det A^I\Delta_I}
\label{3.28}
\eeq
and
\beq
\tau_I(x)=\frac1{\Delta_I} \Wr\left(x^{\gamma_I}\right) \tau_\gamma(x)
\label{9.8}
\eeq
where $\tau_I(t)$ is the tau-function corresponding to the wave function
$\Psi_I(x,z)$ and
$\tau(x)=\tau(x,0,0,\ldots)$.
Applying $a$ to both sides of \eqref{3.28} and using \eqref{1.7} and
\eqref{5.6} we obtain
\beq
\Psi_{aW}(x,z)=\frac{\sum\det A^I\Delta_I
\tau_{a(\gamma+dN\delta_I-n\delta)}(x) \Psi_{a(\gamma+dN\delta_I-n\delta)}(x,z)
} {\sum\det A^I\Delta_I \tau_{a(\gamma+dN\delta_I-n\delta)}(x)}.
\label{9.6}
\eeq
We compute
\beq
a(\gamma +dN\delta_I-n\delta) =a(\gamma) +dN\delta_{I^0} -(dN-n)\delta,
\label{9.7}
\eeq
which is a Darboux transformation of $a(\gamma)$.

The eqs. \eqref{9.8} and \eqref{7.16} imply
\beq
\frac{\tau_I(x)}{\tau_J(x)} =\frac{x^{p_I}}{x^{p_J}}.
\label{9.9}
\eeq
To apply \eqref{9.9} in \eqref{9.6} we have to compute $p_{I^0}$ but for
$a(\gamma)$ instead of $\gamma$. It is a simple exercise to see that
$$
p_{I^0}(a(\gamma))=p_I(\gamma)\equiv p_I.
$$
Using this we obtain
$$
\Psi_{aW}(x,z)=z^{-dN+n} \frac{\sum\det A^I\Delta_I x^{p_I}
L_{(a(\gamma))_{I^0}}} {\sum \det A^I\Delta_I x^{p_I}}
\Psi_{a(\gamma)}(x,z).
$$
Now \prref{4.3} gives (b) because
$$
(L_\beta)^* =(-1)^N L_{a(\beta)} \quad {\rm for}\ \beta\in \Cset^N
$$
and
$$
a\Bigl((a(\gamma))_{I^0}\Bigr)=\gamma_{I^0} -n\delta_{I^0}.
$$
\qed

\smallskip

In the following example we consider the case when there are logarithms in the
basis \eqref{5.20a} of $\ker P$.
\bex{log}
Let $d=2$, $\beta = (1,1,1)$, $\gamma = \beta^2 = (1,1,1,4,4,4)$ and $\ker P$
has a basis
\beqa
&&f_0(x) = x^4,
\nn\\
&&f_1(x) = a_1 x + 2 a_2 x^4 \ln x,
\nn\\
&&f_2(x) = a_0 x + a_1 x \ln x + a_2 x^4 \ln^2 x.
\nn
\eeqa
Using that $\ln^k x = \partial_\epsilon^k x^\epsilon|_{\epsilon=0}$ we
approximate the above functions with
\beqa
&&f_0(x,\epsilon) = x^4,
\nn\\
&&f_1(x,\epsilon) =
a_1 x^{1+\epsilon} + 2 a_2 \epsilon^{-1} (x^{4+\epsilon} - x^4),
\nn\\
&&f_2(x,\epsilon) =
a_0 x^{1+2\epsilon} + a_1 \epsilon^{-1} (x^{1+2\epsilon} - x^{1+\epsilon})
+ a_2 \epsilon^{-2} (x^{4+2\epsilon} - 2x^{4+\epsilon} + x^4).
\nn
\eeqa
Consider the Darboux transformation $W(\epsilon)$ of $V_{\beta(\epsilon)}$,
where $\beta(\epsilon) = (1, 1+\epsilon, 1+2\epsilon)$, with a basis of the
operator $P(\epsilon)$ consisting of the functions $f_k(x,\epsilon)$.
After changing the basis this corresponds to a matrix (cf.\ \eqref{9.5})
$$
A(\epsilon) =
\left(\matrix{
0 & 1 \cr
0 & 0 \cr
0 & 0 \cr
}\right|\left.\matrix{
0 & 0 \cr
\epsilon a_1 & 2a_2 \cr
0 & 0 \cr
}\right|\left.\matrix{
0 & 0 \cr
0 & 0 \cr
\epsilon(\epsilon a_0 + a_1) & a_2 \cr
}\right).
$$
We apply \eqref{3.28} for $\tau_{W(\epsilon)}$. To make the limit
$\epsilon\to0$ we note that the numerator and the denominator depend
polynomially on $\epsilon$ and that (in the notation of \eqref{3.28})
both $\tau_{\{2,3,6\}}$ and $\tau_{\{2,4,5\}}$ tend to one and the same
Bessel tau-function. So after canceling $\epsilon^3$ and setting $\epsilon=0$
we obtain that $\tau_W$ is a linear combination of 3 Bessel tau-functions:
\beq
\tau_W = \frac
{  9a_1^2 \tau_{(-2,1,1,4,4,7)} + 18a_2(a_0 - a_1) \tau_{(-2,-2,1,4,7,7)} +
4a_2^2 \tau_{(-2,-2,-2,7,7,7)} }
{  9a_1^2 + 18a_2(a_0 - a_1) + 4a_2^2 }.
\eeq
As in the proof of \prref{9.1} from this formula one can compute the operators
$P$, $Q$, $P_{\rm b}$ and $Q_{\rm b}$. It is clear that they also can be
obtained by taking the limit $\epsilon\to0$ directly in the corresponding
expressions for $W(\epsilon)$.
\qed
\eex
{}From here to the end of the subsection we shall restrict ourselves to
the case when $\beta^d=\gamma$ has different coordinates. We choose
the following basis of $\ker L_\beta^d$ (cf.\ \cite{MZ})
\beq
\Phi_{(k-1)d+j}(x) := \mu_{kj} x^{\beta_k+(j-1)N}, \quad 1\le k\le N,\
1\le j\le d,
\label{9.10}
\eeq
where
$$\mu_{k,1}:= 1, \quad \mu_{kj} := \mu_{k,j-1}\cdot\prod_{i=1}^N
(\beta_i-\beta_k-(j-1)N)^{-1}.$$
In this basis the action of $L_\beta$ is
quite simple:
\beq
L_\beta\Phi_{(k-1)d+j} =\cases{
\Phi_{(k-1)d+j-1}, &for $2\le j\le d$\cr
0,                 &for $j=1$.\cr}
\label{9.11}
\eeq
Let a basis of $\ker P$ be
\beq
f_k(x)=\sum_{i=1}^{dN} a_{ki}\Phi_i(x),\quad k=0,\ldots,n-1.
\label{9.12}
\eeq
\bex{9.2}
Let $n=d$, $\beta_i-\beta_j\in N\Zset$ for all $i$, $j$ and
the matrix $A=(a_{ki})$ has the form:
\beq
A=\pmatrix{
t_0^{(1)}    &       &       &          &\ldots  & t_0^{(N)}\cr
\noalign{\vskip3pt}
t_1^{(1)}    & t_0^{(1)} &   &          &\ldots  & t_1^{(N)} &t_0^{(N)}\cr
\noalign{\vskip3pt}
t_2^{(1)}    & t_1^{(1)} & t_0^{(1)} &   &\ldots  & t_2^{(N)} & t_1^{(N)}
   & t_0^{(N)}\cr
\vdots  & \vdots &\hfill&\ddots\hfill & &   \vdots & \vdots &\ddots\hfill\cr
\noalign{\vskip3pt}
t_{n-1}^{(1)} & t_{n-2}^{(1)} &\ldots& t_0^{(1)} & \ldots  &
t_{n-1}^{(N)} & t_{n-2}^{(N)} &\ldots& t_0^{(N)}}
\label{9.13}
\eeq
The type of the matrix is tantamount to the identities $L_\beta f_0=0$,
$L_\beta f_{k+1}=f_k$, $k=1,\ldots,n-1$. Then $\ker P$ is invariant under the
action of $L_\beta$ and by \prref{3.5} the operator $L=PL_\beta P^{-1}$ is
differential of order $N$ and solves the bispectral problem. For a generic
$\beta\in \Cset^N$ the spectral algebra has rank $N$ (i.e.\ it is $\Cset
[L]$). This family can be considered as the most direct generalization of the
{\em{``even case''\/}}
of J.~J.~Duistermaat and F.~A.~Gr\"unbaum \cite{DG} (see also \cite{MZ}). When
$N=2$ our example coincides with it but for $N>2$ here we present a completely
new class of bispectral operators.
\qed
\eex
In connection with the above example we prove the following proposition.
\bpr{7.9}
Let $W\in\GrBb$ {\rm(}$\beta\in\Cset^N$--generic{\rm)} be such that $\A_W$
contains an operator of order $N$. Then $W$ is a monomial Darboux
transformation of $V_\beta$, i.e.\ $W\in\GrMBb\cap\Gr^{(N)}$.
\epr
\noindent
\proof
\prref{3.5} implies that $W\in\GrBb$ belongs to $\GRN$ iff
\beq
L_\beta(\ker P)\subset \ker P.
\label{7.29}
\eeq
If we suppose that $W\not\in Gr_{MB}(\beta)$ then $\ker P$ would contain some
elements of the form \eqref{5.20c}. The action of $L_\beta$ on them is easily
computed:
\beqa
L_\beta D^k_\lambda \Psi_\beta(x,\varepsilon^i\lambda)
&=& D^k_\lambda L_\beta\Psi_\beta(x,\varepsilon^i\lambda)
= D^k_\lambda\bigl(\lambda^N\Psi_\beta(x,\varepsilon^i\lambda)\bigr)
\nn\\
&=& \lambda^N(D_\lambda+N)^k \Psi_\beta(x,\varepsilon^i\lambda).
\nn
\eeqa
Thus the linear space $\span\left\{ D^k_\lambda
\Psi_\beta(x,\varepsilon^i\lambda)\right\}_{0\le k\le m}$
can be identified with the space of polynomials in $D$ of degree $\le m$, with
the action of $L_\beta$ corresponding to $P(D)\mapsto \lambda^N P(D+N)$.
It is clear that all the $L_\beta$-invariant subspaces are of the form
$$
\span\left\{ D^k_\lambda
\Psi_\beta(x,\varepsilon^i\lambda)\right\}_{0\le k\le k_0}
$$
for some $k_0$. The corresponding polynomial Darboux transformation is trivial
in the sense that it leads again to the same plane $V_\beta$ (the operator
$P=(L_\beta-\lambda)^{k_0}$ commutes with $L_\beta$).
Therefore $W\in Gr_{MB}(\beta)$.
\qed

\smallskip

In the same manner as in \exref{9.2}, one can build for arbitrary $k$ rank $N$
bispectral algebras with the lowest order of the operators equal to $k N$ .

It is clear that when the matrix $A$ is not of the form \eqref{9.13} (or a
direct sum of such matrices) then $\ker P$
(given by \eqref{9.12})  is not invariant under the action of $L_\beta$.
\prref{3.5} implies that in this case the spectral algebra does not contain
operators of order $N$. The following example is one of the simplest of this
type.
\bex{9.3}
Let $N=2$, $\beta=(\beta_1,\beta_2)$, $\beta_1+\beta_2=1$, $d=n=2$. We take
$\ker P$ with a basis \eqref{9.12} where
$$
A=\left(\matrix{1 & a\cr 0 & 0\cr}\ \right|\left.\matrix{ 0 & 0\cr 1 &
b\cr}\right)
$$
for some $a,b\in \Cset$, i.e.
\begin{eqnarray*}
&&f_0(x) = \Phi_1(x)+a\Phi_2(x)= x^{\beta_1} +\frac{a}{2(\beta_1-\beta_2+2)}
x^{\beta_1+2},\\
&&f_1(x) = \Phi_3(x)+b\Phi_4(x)= x^{\beta_2} +\frac{b}{2(\beta_2-\beta_1+2)}
x^{\beta_2+2}.
\end{eqnarray*}
Then
$$
L_\beta f_0(x)=ax^{\beta_1},\quad
L_\beta f_2(x)=bx^{\beta_2}
$$
and $\ker P$ is not invariant under $L_\beta$ when $ab\not=0$. The spectral
algebra $\A_W= PL^2_\beta\Cset[L_\beta]P^{-1}$ consists of operators of
orders 4, 6, 8, 10, $\ldots$

This example is also interesting for the fact that it does not require
$\beta_1-\beta_2\in 2\Zset$.
The generalization for arbitrary $N$ is obvious.
\qed
\eex
Another example illustrating \prref{3.5} is the following one.
\bex{9.4}
Let $N=2$, $\beta=(\beta_1,\beta_2)\in\Cset^2$, $\beta_1+\beta_2=1$,
$\beta_1-\beta_2\in 2\Zset$, $d=4$, $n=2$.
We take $\ker P$ with a basis \eqref{9.12} where
$$
A=\left(\matrix{\lambda & 0 & 0 & 0\cr 0 & 0 & 1 & 0\cr}\ \right|
\left.\matrix{ \lambda a+b & \lambda b & 0 & 0\cr 0  &0 & a &b\cr}\right)
$$
for some $a,b,\lambda\in \Cset$.

Then it is easy to see that $\ker P$ is invariant under the operator
$L^3_\beta+ \lambda L^2_\beta$ but it is not invariant under any polynomial of
$L_\beta$ of degree $\le2$. On the other hand $\ker P\subset \ker L^4_\beta$
obviously implies $L_\beta^{4+k}\ker P\subset \ker P$ for $k\ge0$. Therefore
the spectral algebra $\A_W$ is the linear span of the operators
$$
P\left(L_\beta^3 +\lambda L_\beta^2\right) P^{-1},\quad
PL_\beta^{4+k} P^{-1}, \quad k\ge0.
$$
This example is interesting for the fact that (for $\lambda\not=0$) the
operator of minimal order in the spectral algebra is not a Darboux
transformation of a power of $L_\beta$, although the Darboux transformation is
monomial.
\qed
\eex
In the last example of this subsection we show that for $d=n=1$ our results
agree with those of \cite{Z}.
\bex{9.5}
Let $d=n=1$, $\ker P=\Cset f_0,$
$$f_0(x)=\sum_{i=1}^N a_i x^{\beta_i}, \qquad
P=\partial_x -\frac{f_0'(x)}{f_0(x)}$$
and $\beta_i-\beta_j\in N\Zset$ if
$a_ia_j\not=0$.
Then
\begin{eqnarray*}
L&=& PL_\beta P^{-1}=PQ=
    \left(\sum a_i x^{p_i}\right)^{-1}
    \left(\sum a_i x^{p_i}\left(\partial_x-\frac{\beta_i}{x}\right)\right)
\times\\
&&\quad\times\left(\sum a_i
\frac{P_\beta(D_x+N)}{D_x+N-\beta_i} x^{p_i-N+1}\right)
    \left(\sum a_i x^{p_i}\right)^{-1},\\
\Lambda&=&P_{\rm b}Q_{\rm b}=
    \left(\sum a_i\left(\partial_z-\frac{\beta_i}{z}\right)(L_\beta)^{p_i/N}
       \right) \times\\
&&\quad\times\left(\sum a_i(L_\beta)^{p_i/N}
x^{-N+1}\frac{P_\beta(D_z+1)}{D_z+1-\beta_i}\right)
\end{eqnarray*}
where $p_i=\beta_i -\beta_{\min}$, $\beta_{\min}=\min\limits_{a_i\not=0}
\beta_i$,
$$P_\beta(D)=\prod _{i=1}^N (D-\beta_i), \qquad D_x=x\partial_x.$$
We have
$$
\Theta(x)=x^N\left(\sum a_i x^{p_i}\right)^2,\quad
\deg \Theta = N+2(\beta_{\max}-\beta_{\min})
$$
where $\beta_{\max}=\max\limits_{a_i\not=0} \beta_i$.
When $f_0(x)= tx^{\beta_1}+x^{\beta_2}$, $\beta_2-\beta_1=N\alpha$, $\alpha\in
\Zset_{\ge0}$
$$
\Theta(x)=x^N(t+x^{N\alpha})^2
$$
and we obtain the operator $\Lambda$ from \cite{Z}.
\qed
\eex
\subsection{Polynomial Darboux transformations}
In this subsection we shall consider the simplest case of polynomial Darboux
transformation of an operator of order $N$, namely when the polynomial $h(z)$
{}from \eqref{3.15} is equal to $(z - \lambda^N)^2$ for some
$\lambda\in\Cset\setminus 0.$ Using the kernels of the operators $P$, $Q^*$,
$P_{\rm b}$ and $Q_{\rm b}^*$ from 
(\ref{3.8}, \ref{3.9}, \ref{6.4}, \ref{6.9}),
we describe the action of the
involutions $a$ and $b$ ($b_1$ in the Airy case). The Propositions \ref{pe1},
\ref{pe2} below raise some interesting questions and conjectures.

The Bessel and Airy cases are very similar. We shall consider first the Airy
one since it is simpler.

Let $\GrA{\alpha}$, $\alpha\in\Cset^{N-1}.$ Set
\beq
h(z) = (z - \lambda^N)^2, \quad g(z) = f(z) = z^N - \lambda^N.
\label{98}
\eeq
Then $\ker h(L_\alpha)$  has a basis of the form
\beq
\Bigl\{ \partial_x^k \Psi_\a(x,\varepsilon^j\lambda)\Bigr\}_
{0\le j\le N-1,\; k=0, 1}
\label{9.14}
\eeq
and  $\ker P$ has a basis
\beq
f_j(x) =
\Psi_\a(x,\varepsilon^j\lambda) + a \partial_x
\Psi_\a(x,\varepsilon^j\lambda),
\quad 0\le j\le N-1
\label{99}
\eeq
for some $a\in\Cset$.

We shall start with the case $N=2$.
The following example is due to \cite{G3, LP}. We shall obtain it
as the simplest special case of \thref{8.5}.
\bex{9.9}
Let $N=2$ and $\a=(\a_0)\in\Cset^1$.
For fixed $\a_0, a\in \Cset\setminus 0$ we take the basis \eqref{99}
of $\ker P$:
\beq
f_k(x) =  \psi_k(x) + a\partial_x \psi_k(x),
\quad k=0,1
\label{100}
\eeq
where $\psi_k(x) = \Psi_\a(x,(-1)^k \lambda)$.
Using that
\beqa
&& \partial_x f_k =  a(\a_0 x + \lambda^2)\psi_k + \partial_x \psi_k,
\nn\\
&& \partial_x^2 f_k =
  (a\a_0 + \a_0 x + \lambda^2)\psi_k + a(\a_0 x + \lambda^2)\partial_x \psi_k
\nn
\eeqa
we compute $P$ from
\beqa
P\varphi &=& \frac{\Wr(f_0,f_1,\varphi)}{\Wr(f_0,f_1)}
\nn\\
&=& \frac{\left|\matrix{
  1 & a & \varphi\cr
  a(\a_0 x + \lambda^2) & 1 & \partial_x\varphi\cr
  a\a_0 + \a_0 x + \lambda^2 & a(\a_0 x + \lambda^2) & \partial_x^2\varphi
}\right|}
{\left|\matrix{
  1 & a \cr
  a(\a_0 x + \lambda^2) & 1
}\right|}.
\nn
\eeqa
The result is
\beq
P=\partial_x^2 + \frac{a^2\a_0}{1-a^2(\a_0 x + \lambda^2)} \partial_x
+\frac{a^2(\a_0 x + \lambda^2)^2
- (\a_0 x + \lambda^2) - a\a_0}{1-a^2(\a_0 x + \lambda^2)}.
\label{101}
\eeq
This expression coincides with that given in \cite{G3} if we set
$$
\a_0 = \frac2{2+3t}=\frac2s,\quad a=\frac s{2y},\quad \lambda=0.
$$
We compute the operators $P$, $Q$ and $Q^*$ as follows.
If we write
\beqa
&&P = \partial_x^2 + p_1(x) \partial_x + p_0(x)
\nn\\
&&Q = \partial_x^2 + q_1(x) \partial_x + q_0(x)
\nn\\
&&Q^* = \partial_x^2 + \widetilde q_1(x) \partial_x + \widetilde q_0(x)
\nn
\eeqa
then the identity $Q P = h(L_\a)$ imply
$$
q_1 + p_1 = 0, \;\; 2p_1' + q_1  p_1 + p_0 + q_0 = -2(\a_0 x + \lambda^2)
$$
and
$$
\widetilde q_1 = - q_1, \quad \widetilde q_0 = -q_1' + q_0.
$$
Our observation is that because $P^* Q^* = h(L_{a(\a)})$ and $\Psi_{a W} =
f^{-1} Q^* \Psi_{a(\a)}$, the operator $Q^*$ has a basis
of the form \eqref{100} with some $b\in\Cset$ instead of $a$ and $a(\a)$
instead
of $\a$. Comparing the above expressions for $Q^*$ with \eqref{101} we obtain
that $b = -a$.
By \thref{8.5} the operator $P_{\rm b}$ also has a basis \eqref{9.14} with
some $c$ instead of $a$ and $\mu$ instead of $\lambda$. On the other hand we
can compute it directly using eqs.\ (\ref{8.15}, \ref{8.16}).
Then $g_{\rm b}(z) = 1 - a^2 (z^2 + \lambda^2)$ which on the other hand is
up to a constant $z^2 - \mu^2$. This gives
\beq
\mu^2 = \frac{1-a^2\lambda^2}{a^2}.
\label{102}
\eeq
The other coefficients give a surprising result: $c=a$.
In conclusion, if we denote the operator $P$ from \eqref{101} with
$P(a,\lambda)$ then
\beq
P = P(a,\lambda), \quad Q = P^*(-a,\lambda), \quad
P_{\rm b} = P(a,\mu), \quad Q_{\rm b} = P^*(-a,\mu)
\label{103}
\eeq
where $\mu$ and $\lambda$ are connected by \eqref{102}.
\qed
\eex
The next example is completely analogous to the above one but to the best of
our knowledge it is new.
\bex{9.10}
For $N=3$ the Airy operator is
$$
L_\alpha=\partial_x^3+\alpha_2 \partial_x -\alpha_0 x,
$$
$\alpha=(\alpha_0,\alpha_2)\in \Cset^2$, $\alpha_0\not=0$.
We take $P$ with a basis \eqref{99} ($N=3$).
Then using the eq.\ \eqref{8.8}  we compute
\begin{eqnarray*}
P &=& \partial_x^3 -\frac{a^3\a_0}{a^3(\a_0 x+\lambda^3)
+(1+a^2\a_2)}\,\partial_x^2\\
&+& \frac{a^3\a_2(\a_0 x+\lambda^3)+(1+a^2\a_2)\a_2 +a^2\a_0}
      {a^3(\a_0 x+\lambda^3)+(1+a^2\a_2)}\,\partial_x\\
&-& \frac{a^3(\a_0 x+\lambda^3)^2 + a\a_0(1+a^2\alpha_2)
(1+a^2\alpha_2)(\a_0 x+\lambda^3)}
       {a^3(\a_0 x+\lambda^3)+(1+a^2\a_2)}.
\end{eqnarray*}
A direct computation using \prref{4.3}, \thref{8.5} and $Q P = h(L_\a)$ leads
to
\beq
P = P(a,\lambda), \quad Q = -P^*(-a,-\lambda), \quad
P_{\rm b} = P(a,\mu), \quad Q_{\rm b} = -P^*(-a,-\mu)
\label{103'5}
\eeq
with $\mu$ given by
\beq
\mu^3 + \lambda^3 = -\frac{1+a^2 \a_2}{a^3}
\eeq
\qed
\eex
The above examples can be generalized for arbitrary $N$ as follows.
\bpr{e1}
Denote by $P=P(a,\lambda)$ the operator $P$ with a basis \eqref{99}.
Then in the above notation we have
\beq
Q = (-1)^N P^*(-a,-\lambda), \quad
P_{\rm b} = P(a,\mu), \quad Q_{\rm b} = (-1)^N P^*(-a,-\mu)
\label{103'6}
\eeq
with $\lambda$ and $\mu$
connected by
\beq
\lambda^N + \mu^N = P_{\a'}(-1/a)
\label{104}
\eeq
where $P_{\a'}$ is the polynomial from \eqref{8.1}.
The spectral algebras
\beqa
&&\A_W = P\left(L_\a-\lambda^N\right)^2 \Cset[L_\a] P^{-1}
\label{105}\\
&&\A_{b_1 W} = P_{\rm b}\left(L_\a-\mu^N\right)^2 \Cset[L_\a] P^{-1}_{\rm b}
\label{106}
\eeqa
consist of operators of orders $2N,3N,4N,\ldots$
\epr
\proof
Because
\beqa
&&(-1)^N P^* \, (-1)^N Q^* = (L_{a(\a)} - (-\lambda)^N)^2,
\quad
\Psi_{a W}(x,z) = \frac{(-1)^N Q^* \Psi_{a (\a)}(x,z)}{z^N - (-\lambda)^N},
\nn\\
&&Q_{\rm b} \, P_{\rm b} = (L_\a - \mu^N)^2,
\quad
\Psi_{b W}(x,z) = \frac{P_{\rm b} \Psi_\a(x,z)}{z^N - \mu^N}
\nn
\eeqa
we see that $Q = (-1)^N P^*(b,-\lambda)$ and $P_{\rm b} = P(c,\mu)$
for some $b,c,\mu$.
Using the eq.\
$
L_\a(x,\partial_x)\Psi_\a(x,\varepsilon^j\lambda)
= \lambda^N\Psi_\a(x,\varepsilon^j\lambda)
$
we compute $p_N(x) = \Wr(f_0,f_1,\ldots,f_{N-1})$ as in the proof of
\leref{5.9}. We obtain
$$
p_N(x) =  - (-a)^N (\a_0 x + \lambda^N - P_{\a'} (-1/a)).
$$
Eq.\ \eqref{8.16'5} leads to \eqref{104} because
$g_{\rm b}(z) = \const\cdot (z^N-\mu^N).$
Applying \eqref{104} for $P_{\rm b}$ instead of $P$, we obtain
$$
P_{\a'}(-1/a) = P_{\a'}(-1/c).
$$
Note that the map $a \mapsto c$ is an automorphism of $\Cset{\mathbb P}^1$
since it is an involution.
The only solution of the above equation with this property is $c=a$.

To compute $(-1)^N Q^*$,
 we note that its second coefficient is equal to that of
$P$ which is equal to
$
-p_N'(x)/p_N(x).
$
This, \prref{4.3} and \leref{8.7} imply
$$
\frac{\a_0}{\a_0 x + \lambda^N - P_{\a'} (-1/a)}
= \frac{a(\a)_0}{a(\a)_0 x + (-\lambda)^N - P_{a(\a)'} (-1/b)}
$$
which leads to a polynomial equation for $b$ in terms of $a$ and $\a.$
Because $a \mapsto b$ is an automorphism of $\Cset{\mathbb P}^1$
we obtain that $b=-a$.

The eqs.\ (\ref{105}, \ref{106}) follow from  \prref{7.9}.
\qed

\smallskip\noindent

We shall find the analog of \prref{e1} in the Bessel case. We use the notation
{}from the beginning of the subsection with $\beta\in\Cset^N$ instead of $\a$
and eq.\ \eqref{99} modified as follows (cf. \eqref{5.20c})
\beq
f_j(x) =
\Psi_\beta(x,\varepsilon^j\lambda) + a D_x
\Psi_\beta(x,\varepsilon^j\lambda)
\label{107}
\eeq
($j=0,\ldots,N-1,\; D_x=x\partial_x$).

In the next example we shall study the simplest case $N=2$.
\bex{9.8}
For $N=2$, $\beta=(1-\nu,\nu)$ the corresponding Bessel
operator is
$$
L_\beta=x^{-2}(D_x-(1-\nu))(D_x-\nu)= \partial_x^2 +\frac{\nu(1-\nu)}{x^2},
\quad D_x=x\partial_x.
$$
Using \eqref{2.17'2}
we compute the operator $P$ from
$$
P\varphi =\frac1{x^2}
\frac{\left|\matrix{ f_0 & f_1 &\varphi\cr
                     D_x f_0 & D_xf_1 & D_x\varphi\cr
               D^2_x f_0 & D^2_x f_1 & D^2_x \varphi}\right|}
{\left|\matrix{ f_0  & f_1\cr
               D_xf_0 & D_x f_1}\right|}.
$$
The answer is the following. If we set
$$
\mu^2=\frac{a+1+a^2\nu(1-\nu)}{a^2\lambda^2}
$$
then
$$
P=\frac1{x^2 p_2(x^2)}
\Bigl\{ p_2(x^2) D^2_x + p_1(x^2) D_x + p_0(x^2)\Bigr\}
$$
with
$p_2(x^2)=x^2-\mu^2$, $p_1(x^2)=\mu^2-3x^2$
and
$
p_0(x^2)= -\lambda^2 x^4 + (2\lambda^2\mu^2 + (a+1)(2a-1)a^{-2})x^2
 + ((a+1)a^{-2} - \lambda^2\mu^2)\mu^2.
$
The operator $P_{\rm b}$ is (cf.\ \eqref{6.6})
$$
P_{\rm b}=\frac1{g(x)} \Bigl\{ D_x^2 p_2(L_\beta) +D_x p_1(L_\beta)
+p_0(L_\beta)\Bigr\}
$$
and $g_{\rm b}(z) = z^2(z^2 - \mu^2)$.
A straightforward computation shows that if we set $P=P(a,\lambda,\mu)$ then
\beq
Q = P^*(-a/(a+1),\lambda,\mu), \;\;
P_{\rm b} = P(a,\mu,\lambda) L_\beta, \;\;
Q_{\rm b} = L_\beta P^*(-a/(a+1),\mu,\lambda).
\label{108}
\eeq
Therefore we can take
\beq
P_{\rm b} = P(a,\mu,\lambda), \quad
Q_{\rm b} = P^*(-a/(a+1),\mu,\lambda)
\label{109}
\eeq
i.e.\ the involution $b$ acts simply by exchanging $\lambda$ with $\mu$ and
vice versa, while the involution $a$ acts as $a \mapsto -a/(a+1)$.
\qed
\eex
The action of the involutions for arbitrary $N$ is given in the next
proposition.
\bpr{e2}
Denote by $P=P(a,\lambda)$ the operator $P$ with a basis \eqref{107}.
Then we can take $P_{\rm b}$ and $Q_{\rm b}$ such that
\beq
Q = (-1)^N P^* (b,-\lambda),
\quad
P_{\rm b} = P(a,\mu),
\quad
Q_{\rm b} = (-1)^N P^* (b,-\mu)
\label{110}
\eeq
with $\lambda$, $\mu$  and $a$, $b$ connected by
\beq
\lambda^N  \mu^N = P_{\beta}(-1/a),
\quad
\frac{1}{a} + \frac{1}{b} + N-1 = 0
\label{111}
\eeq
where $P_{\beta}$ is the polynomial from \eqref{2.13}.
The spectral algebras
\beqa
&&\A_W = P\left(L_\beta-\lambda^N\right)^2 \Cset[L_\beta] P^{-1}
\label{112}\\
&&\A_{bW} = P_{\rm b}\left(L_\beta-\mu^N\right)^2 \Cset[L_\beta] P^{-1}_{\rm
b}
\label{113}
\eeqa
consist of operators of orders $2N,3N,4N,\ldots$
\epr
\proof
We have
$g_{\rm b}(z) = \const \cdot z^N (z^N - \mu^N)$
for some $\mu\in\Cset$.
Using \eqref{6.7} we compute
$$
g_{\rm b}(z) = z^N \det (D_z^i f_j(z))_{i,j=0,\ldots,N-1}
= (-a)^N z^N (z^N \lambda^N - P_\beta(-1/a))
$$
which gives the value of $\mu$.
We have to prove that $P_{\rm b}$ given by \eqref{6.6} (which is of order $2N$)
is divisible by $L_\beta$ from the right. Indeed, it is easy to see that
$$
P_{\rm b}(x,\partial_x) x^{\beta_i} =
P(x,\partial_x) x^{\beta_i}|_{\lambda=0}
$$
and the proof of \leref{5.9} implies $P|_{\lambda=0} = L_\beta$.
Thus we can take $P_{\rm b} = P(c,\mu)$ for some $c\in\Cset$.
Now \eqref{111} implies $P_\beta(-1/a) = P_\beta(-1/c)$ leading to $c=a$.

Finally, as in the Airy case, if $Q = (-1)^N P^*(b,-\lambda)$ for some
$b\in\Cset$ then
$$
P_\beta(-1/a) = (-1)^N P_{a(\beta)}(-1/b) = P_\beta(1/b + N-1)
$$
showing that $a^{-1} + b^{-1} + N-1 = 0$.

The eqs.\ (\ref{112}, \ref{113}) follow from \prref{7.9}.
\qed

\smallskip

In conclusion we want to make some comments.
In the case $N=1$ the adjoint involution $a$ has a simple and beautiful
geometric interpretation (see \cite{W}): in terms of Krichever's
construction it preserves the spectral curve and maps the ``sheaf of
eigenfunctions'' into some kind of a dual sheaf. In \cite{W} G.~Wilson also
posed the problem of describing the action of the bispectral involution on
$\Gr^{ad}$.

We think that in the general case the study of the action of the involutions
$a$ and $b$ on the bispectral manifolds of polynomial Darboux transformations
of Bessel and Airy planes is equally interesting and difficult task.

The above examples lead us to the conjecture that the involutions $a$ and $b$
($b_1$ in the Airy case) possess some universality property.
Any polynomial Darboux transformation $W$ of a Bessel plane $V_\beta$
(respectively an Airy plane $V_\alpha$) is determined by the points
$\lambda_1,\ldots,\lambda_N$ ($\not= 0$) at which the conditions $C$ are
supported (see \eqref{5.1}), by the matrix $A$ defined by \eqref{5.20c}
(resp.\ \eqref{5.20b}), and of course by the vector $\beta$ (resp.\
$\alpha$). Then the corresponding matrices for $a W$ and $b W$ (resp.\ $b_1 W$)
depend only on the matrix $A$. The point is that they do not depend on the
points
$\lambda_1,\ldots,\lambda_N$ at which the conditions $C$ are
supported nor on the vector $\beta$ (resp.\ $\alpha$).
\begin{small}
\renewcommand{\refname}{ {\flushleft\normalsize\bf{References}} }
    
\end{small}
\end{document}